\documentclass[%
aip,cha,author-year,
 amsmath,amssymb,
 reprint,%
 single-column,
]{revtex4-1}

\usepackage{graphicx}        
  \def\thecomma{\ifx,\thenext \else\ifx;\thenext \else\ifx.\thenext \else\ifx!\thenext \else\ifx:\thenext \else \  \fi\fi\fi\fi\fi}
\def\condblank{\futurelet\thenext\thecomma}
\def\ie{{\it i.e.},\condblank}
\def\eg{{\it e.g.},\condblank}
\def\CC#1{\begin{center}{#1}\end{center}}
\def\xx#1{{#1}}
\def\fref#1{Fig.~{\ref{#1}}}
 \def\HH{\textbf{H}}
 \def\VV{\textbf{V}}
\def\ssection#1{\section{#1}}

\usepackage{color}

\usepackage[utf8]{inputenc}
\usepackage[T1]{fontenc}
\usepackage{mathptmx}
\usepackage{etoolbox}
\makeatletter
\def\@email#1#2{%
 \endgroup
 \patchcmd{\titleblock@produce}
  {\frontmatter@RRAPformat}
  {\frontmatter@RRAPformat{\produce@RRAP{*#1\href{mailto:#2}{#2}}}\frontmatter@RRAPformat}
  {}{}
}%
\makeatother

\makeatother
\begin{document}

\title[Reflections on a tube]{{An Introduction to the Unpublished Book\\``Reflections on a
   Tube'' by Mitchell J. Feigenbaum}}
\author{Jean-Pierre Eckmann}
\affiliation{D\'epartement de Physique Th\'eorique and\\Section de
  Math\'ematiques, University of Geneva, Geneva,
  Switzerland}
\date{\today}
\begin{abstract}
  This paper is an adaptation of the introduction to a book
  project by
the late Mitchell J. Feigenbaum (1944-2019). While Feigenbaum is
certainly mostly known for his theory of period doubling cascades, he
had a lifelong interest in optics. His book project is an extremely
original discussion of the apparently very simple study of anamorphs,
that is, the reflections of images on a cylindrical mirror. He
observed that there are \emph{two images} to be seen in the tube, and
discovered that the brain preferentially chooses one of them. I edited
and wrote an introduction to this planned  book. As the
book is still not published, I have now adapted my introduction as a standalone
article, so that some of Feigenbaum's remarkable work will be accessible
to a larger audience.
\end{abstract}

\maketitle

{\bf 
 The late Mitchell J. Feigenbaum (1944-2019) left us with an
 unfinished book whose title is ``Reflections on a Tube.'' While Feigenbaum is
certainly mostly known for his theory of period doubling cascades, he
had a lifelong interest in optics. In the book, he starts with the
study of the image you can see in a vertically placed cylindrical mirror,
usually known as anamorph.
He
observed that there are \emph{two images} to be seen in the tube, and
discovered that the brain preferentially chooses one of them.
Fanning out from this observation, he touches on several associated
problems: What fish see from under the water, the quality of the fish
eye versus the land-animal eye, and many others. As the
book is still not published, I have now adapted my introduction the
book as a standalone
article, so that some of Feigenbaum's remarkable work will be accessible
to a larger audience.
}

\ssection{The book}When Mitchell J. Feigenbaum passed away in 2019,
he left a manuscript of his work on ``Reflections on a tube'' to me,
and, in different
versions, to some of his other friends.
Since everyone in this group knew that Mitchell and I had discussed
many times over the years all aspects of the book, the general feeling
was that I 
should finish the manuscript and publish it as a book under Feigenbaum's name
with myself listed as the editor who also completed the missing pieces. 
The present manuscript is a modified version of my planned introduction to that book. Unfortunately, the book project had to be
given up for difficulties with the copyright which stays with the
heirs. I therefore decided to at least make my introduction available
to others. I still hope that the book will finally appear in
some form or other, but in the meantime, I hope that my
``introduction'' will make Feigenbaum's ideas known to a larger public.
In the meantime there are two papers available in which Gemunu Gunaratne
and I tried to explain some details of Feigenbaum's
work \citep{ruler,witchball}, so his ideas can be followed on a more
technical level.



\makeatletter
\renewcommand\thefigure{\@arabic\c@figure}
\makeatother

\ssection{The subject of the book}

\begin{figure}[ht!]
  \CC{\includegraphics[width=\columnwidth]{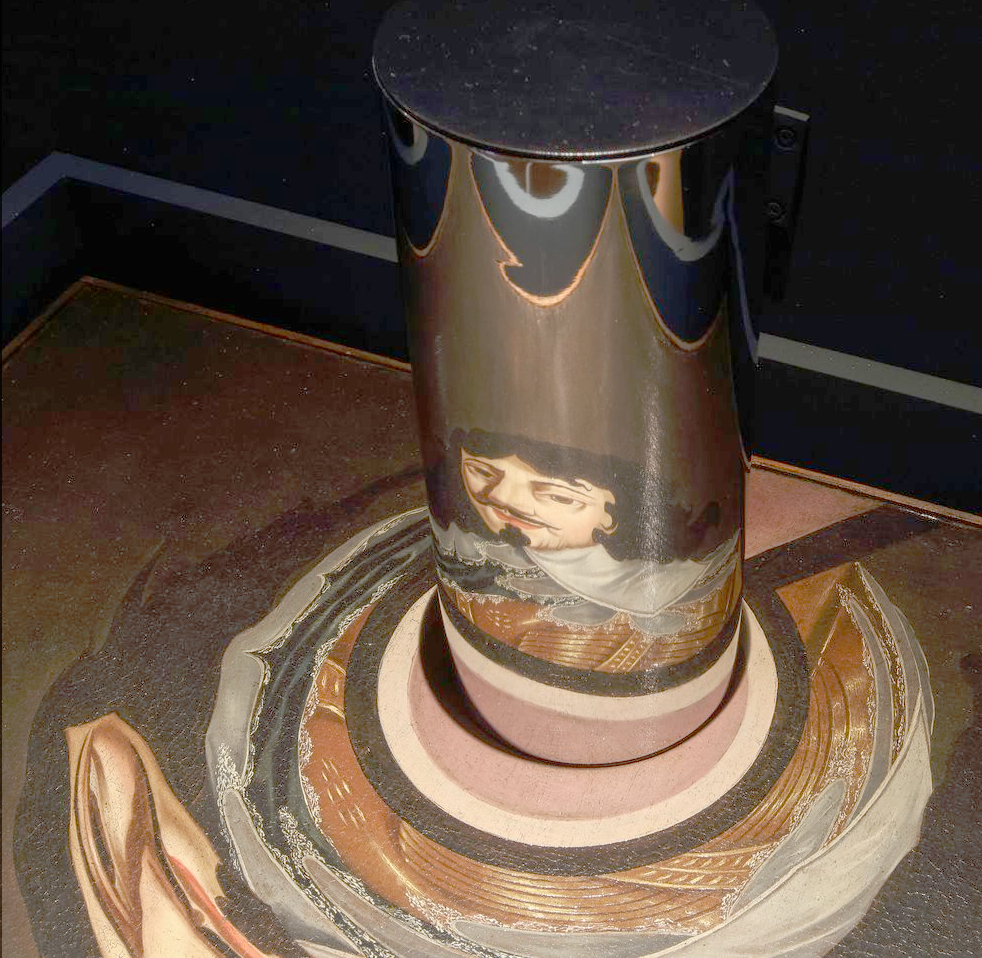}}
    \caption{An anamorph by Jean-Fran\c{c}ois \xx{Niceron} (1616-1646). It
      shows king Louis XIII. (Photo: Alberto Novelli (contrast enhanced).)}\label{fig:anamorph}
    \end{figure}

\begin{figure}[ht!]
  \CC{\includegraphics[width=\columnwidth]{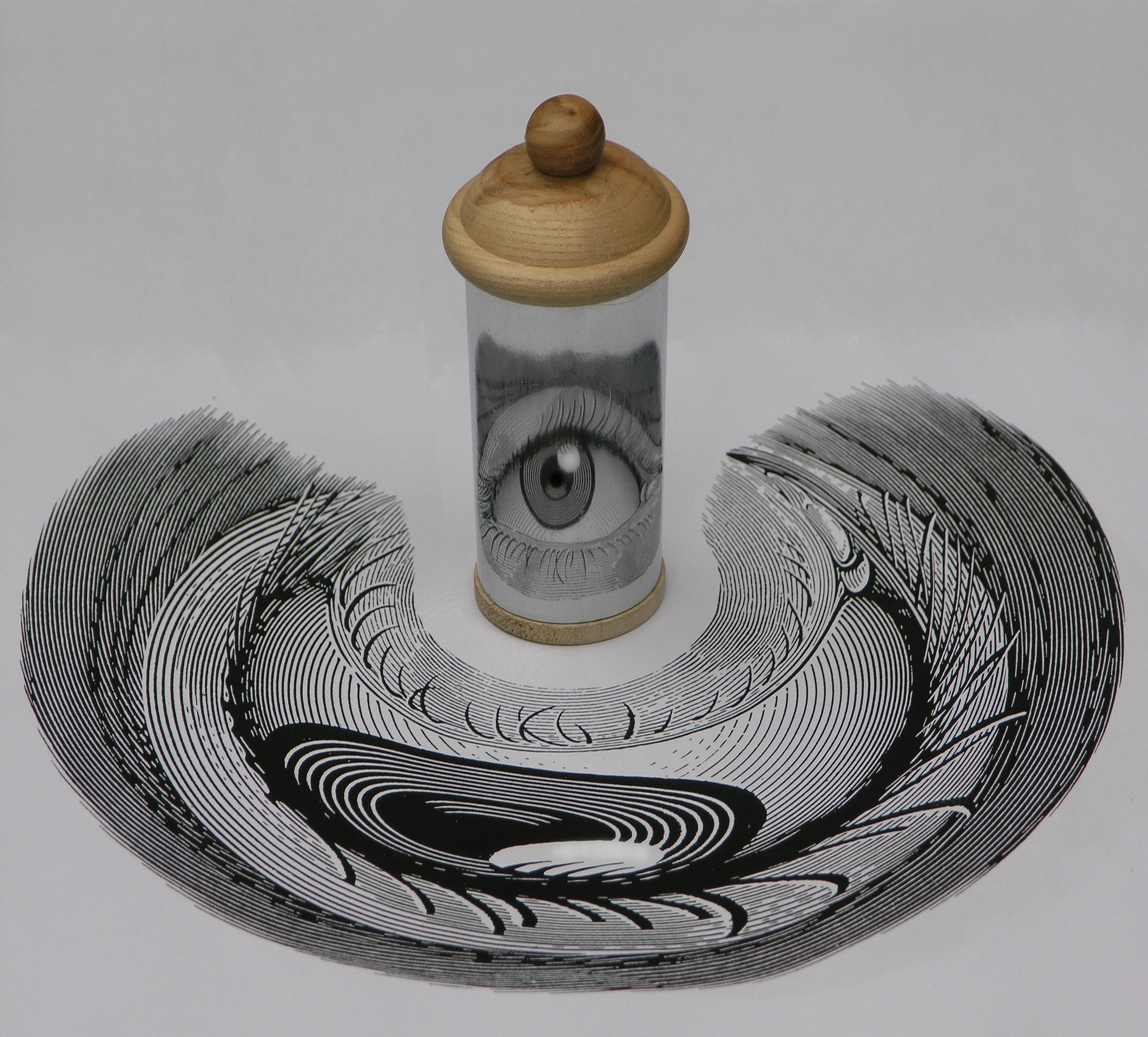}}
  \caption{An anamorph by Istv\'an  Orosz}\label{fig:istvan}
  
\end{figure}

This book is about \emph{anamorphs}\index{anamorph}, reflections of images in a
cylindrical tube. They are known to a large public, from first historical
examples, such as \fref{fig:anamorph} \citep{niceron1638} to modern
works of art, 
such as \fref{fig:istvan}. A drawing is deformed in such a way on a
piece of paper so that the observer will see the undistorted image,
when looking at the tube. In \fref{fig:anamorph} one sees Louis XIII, and
in \fref{fig:istvan} a beautiful eye appears.
Given the many anamorphs one can find, one feels that
their theory must have been
extremely well-studied in optics.
Indeed,
one can find
many programs which allow one to generate the anamorphic picture
on the ground from any sample image.
The novelty of Feigenbaum's work is that, upon studying the visual
properties of the reflections 
in detail, one finds that the theory of anamorphs requires concepts that go way beyond such a
seemingly simple toy problem.
Specifically, Feigenbaum worked out an intriguing
dichotomy of possible interpretations of what one can see. This
dichotomy gives us a glimpse into the inner workings of the human
visual system and its connection to the brain. Feigenbaum's
observations are largely unexplained from a physiological point of view. What I like
about this work is its methodology, which shows how a careful
calculation (in this case in optics) can lead to unexpected observations in
another field (in this case, perception).

\begin{figure}
  \CC{\includegraphics[width=0.6\columnwidth]{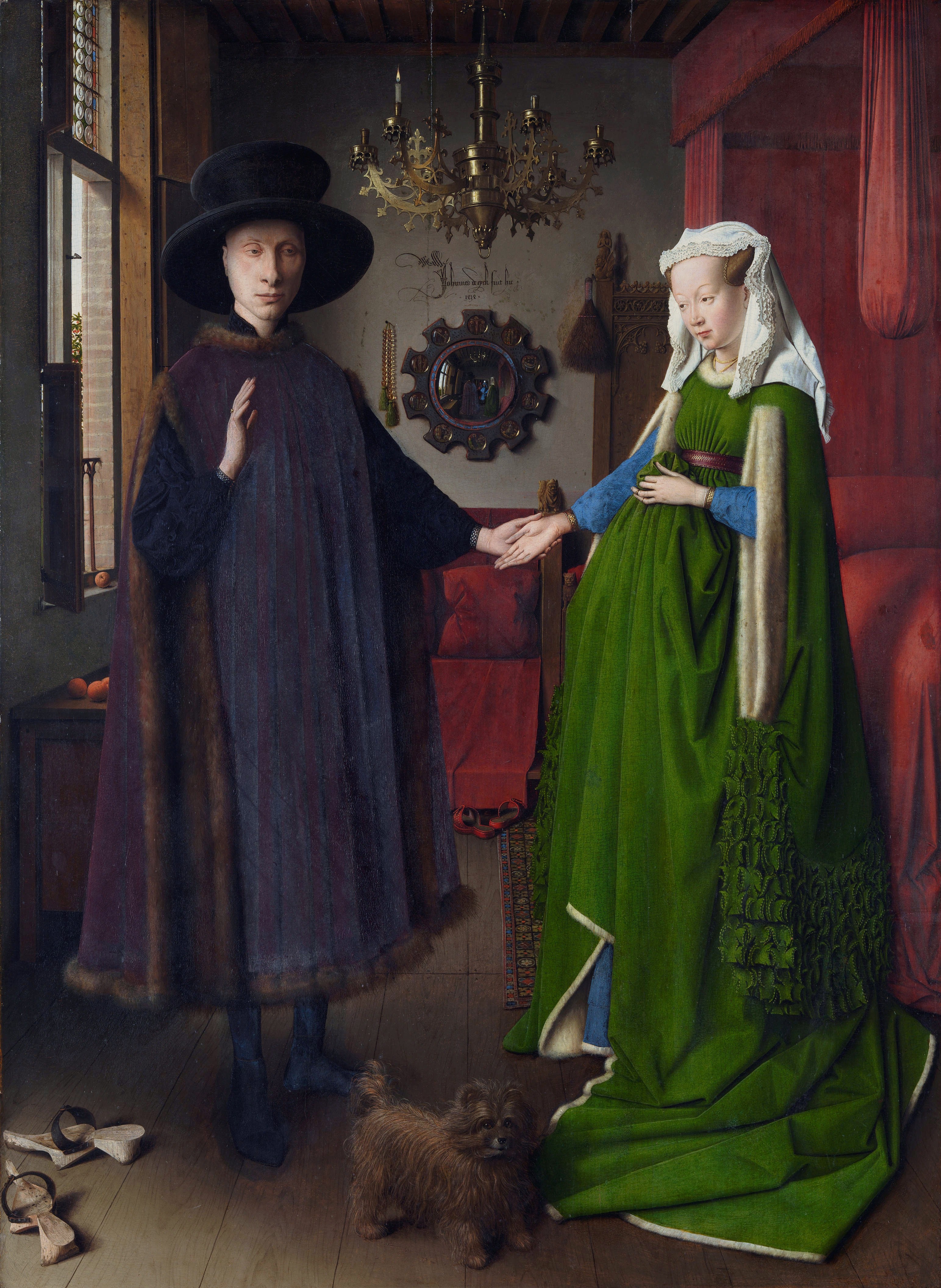}}
   \CC{\includegraphics[width=0.5\columnwidth]{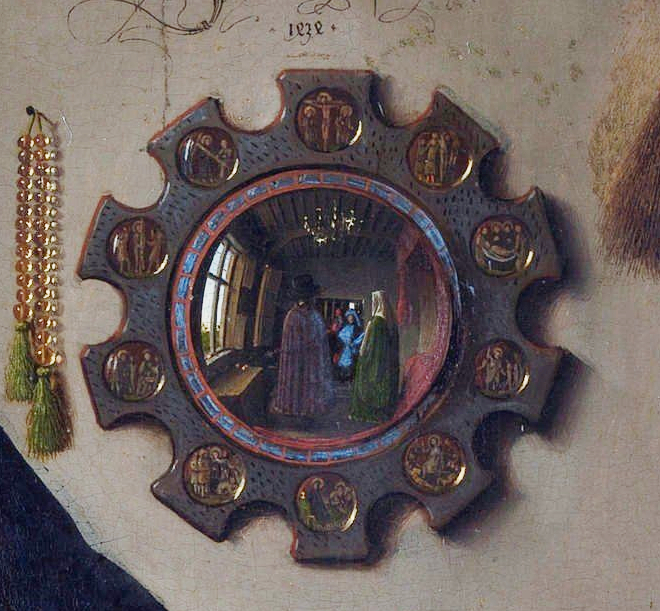}}
  \caption{The Arnolfini portrait of van Eyck 1434 (National Gallery
    London) is considered the first painting with reflections from a
    non-flat mirror. Both views were artificially made lighter for
    better visibility. (Source: Wikimedia commons.)}\label{fig:london}
  
\end{figure}

Mirrors come in many forms: The standard mirror on a wall is flat, but
a mirror can also
be bent  like
the cylinder, rippled like the surface of water, or willfully distorted
 as in \fref{fig:london}. Still, all these  mirrors are two dimensional surfaces\label{pag:2d}. The
theory of his book covers visual aspects of such mirrors.

\ssection{Three possible anamorphs}

The mathematical finding of Feigenbaum is that there are really 3
possible images to be seen in the tube:  namely, a standard one, which he calls ``erect,'' and
two others, which he calls ``3D'' and ``flat.''
\begin{figure}
   {\includegraphics[width=\columnwidth]{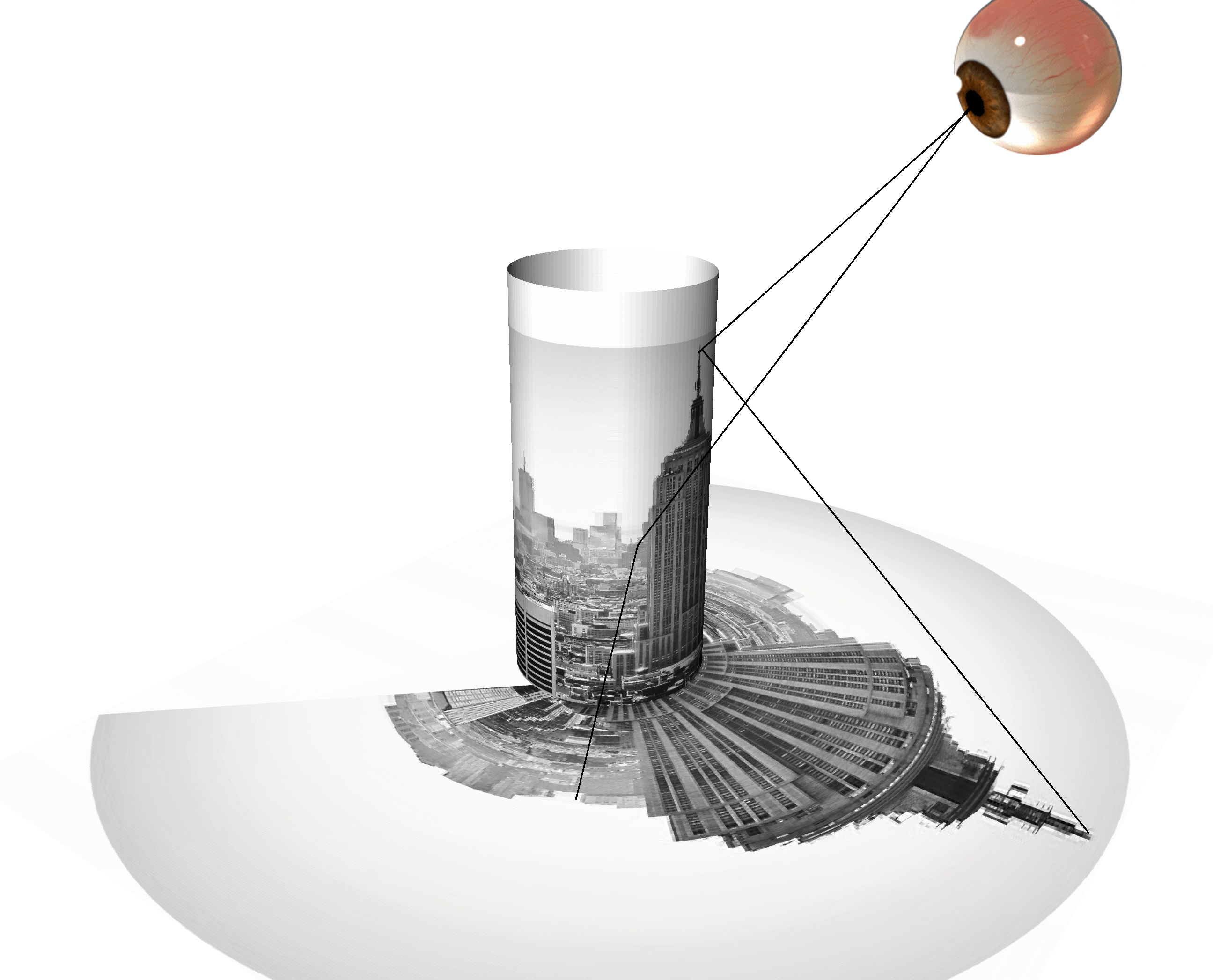}} 
  \caption{The ``erect'' anamorph, where the construction is computed
    as if the image were wrapped on the cylinder.}\label{fig:erect}
  
\end{figure}

The erect one is shown in \fref{fig:erect}. One wraps an image
around the cylinder, fixes the position of the eye, and then draws lines
from the eye to the cylinder and then to the table, using the rules of
reflection. This is what is done in the anamorphs of
\fref{fig:anamorph} and \fref{fig:istvan}. While this procedure leads to
appealing anamorphs, it is actually not correct, because, as is
shown in Feigenbaum's book, these views have no power in the sense that they would be what is
seen by a pin-hole camera, but not by the human eye, which has a
non-negligible opening of the pupil.
\begin{figure}[ht!]
  \CC{\includegraphics[width=0.45\columnwidth]{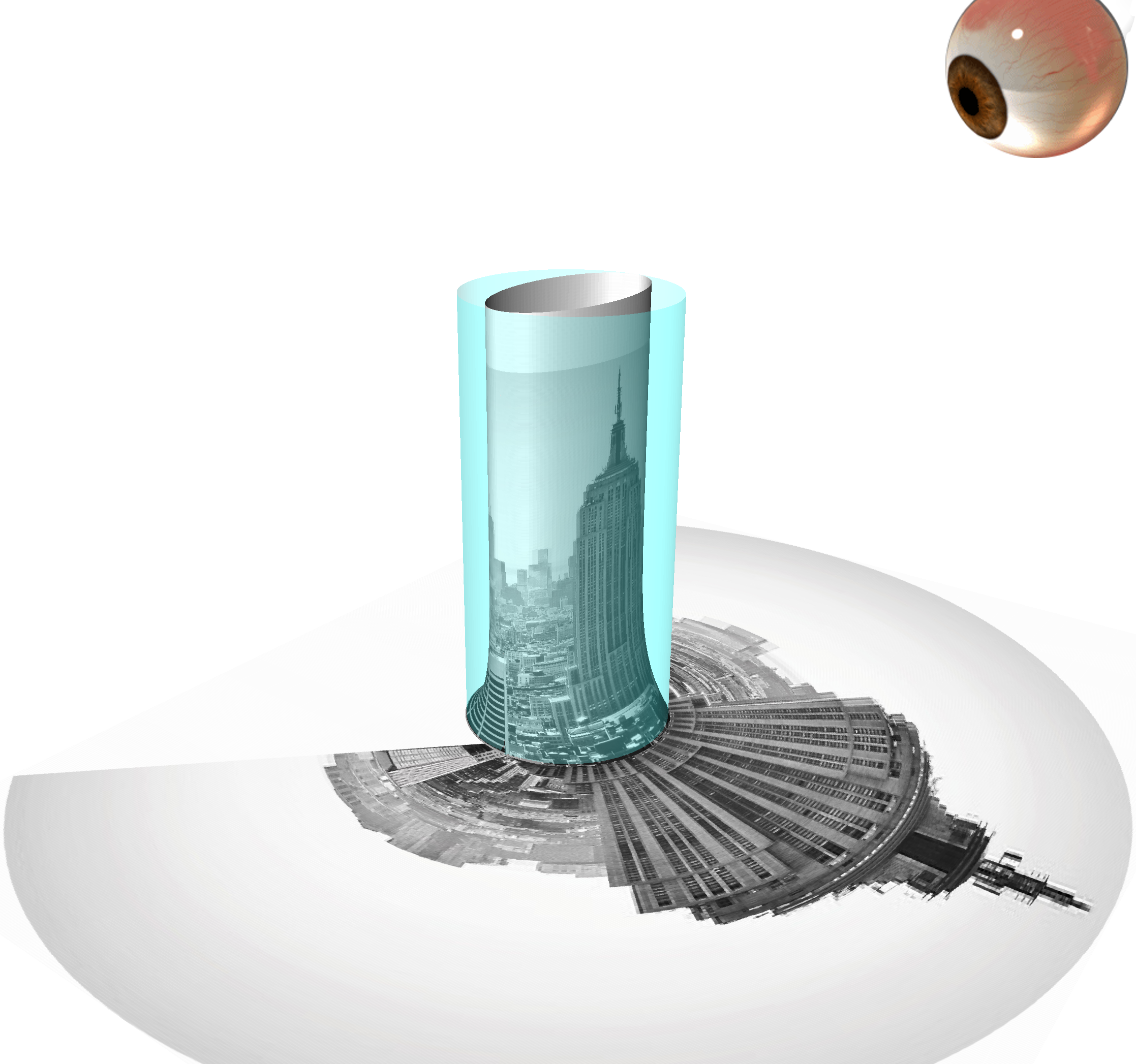}
  \includegraphics[width=0.45\columnwidth]{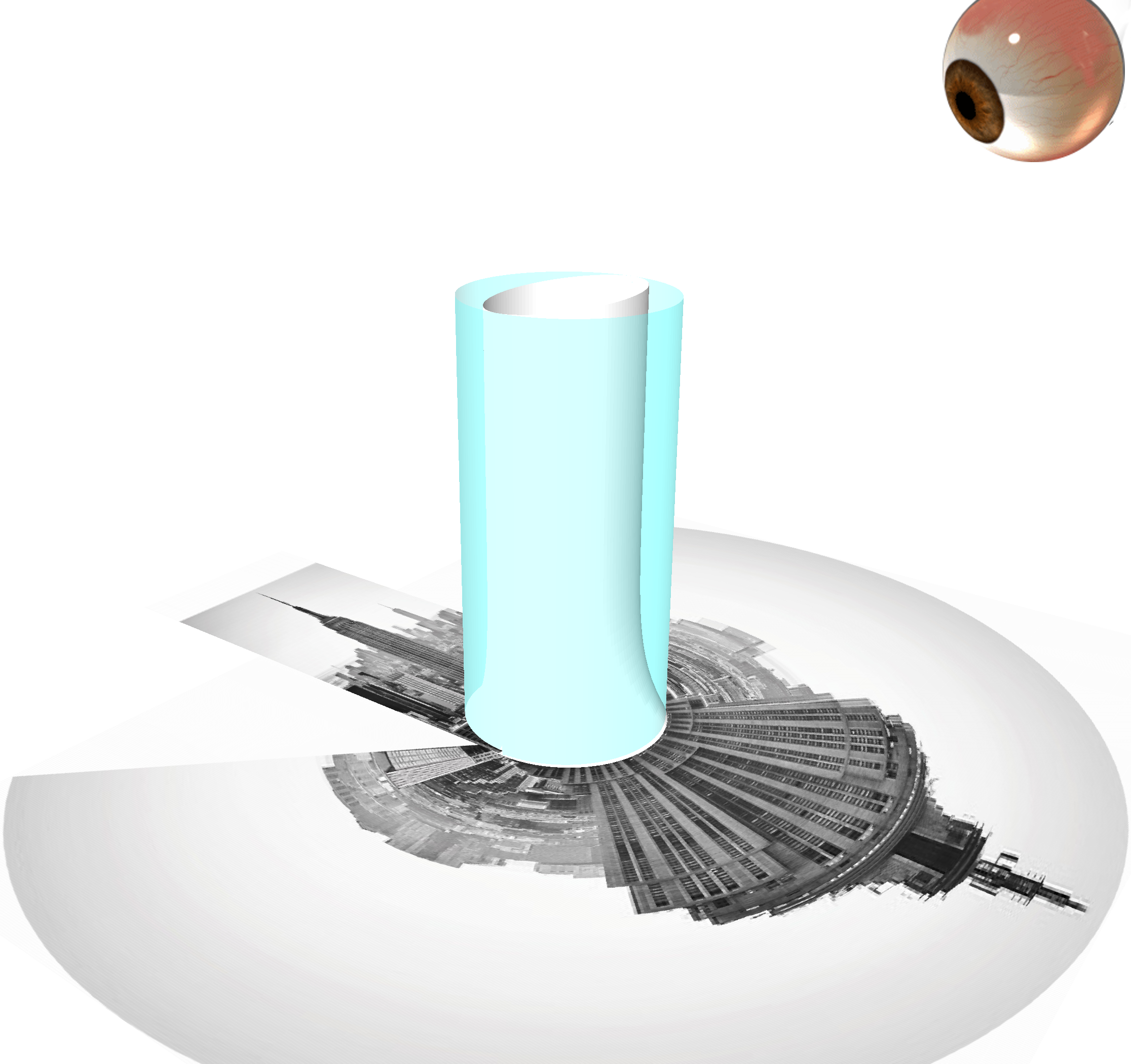}}
  \caption{The 3D and flat anamorphs: Inside the cylinder (left) and flat
    on the table (right). While, of course, the light rays always get
    reflected exactly as shown in \fref{fig:erect}, the virtual image will
    appear not \emph{on}  the surface of the cylinder, but either on a
    surface \emph{inside} the
    cylinder or flat on  the table \emph{behind }it. This is what the observer will really see.}\label{fig:correct}
  
\end{figure}

Using the eye, and not a pinhole camera, one can actually see two
different images, as suggestively shown in
\fref{fig:correct}. The first image
appears on a surface which is in the \emph{interior} of the tube,
while the second  lies \emph{flat} on the table.
Thus, \emph{two} different views are presented to the eye.

A poor man's explanation for the two images is understood, indirectly,
because a 2-dimensional surface has, in every
point, {\emph{two main curvatures}}\footnote{The flat mirror is
  exceptional in this respect, since all directions have the same curvature.}.
 For example, the tube is flat in the
vertical direction and maximally curved in the horizontal direction.

Using the theory of caustics, which will be illustrated below,
Feigenbaum showed that there are indeed two images, as in \fref{fig:correct},
both of which have more intensity  than the erect image sketched in
\fref{fig:erect}.
\begin{figure}[h!]
  \CC{\includegraphics[width=\columnwidth]{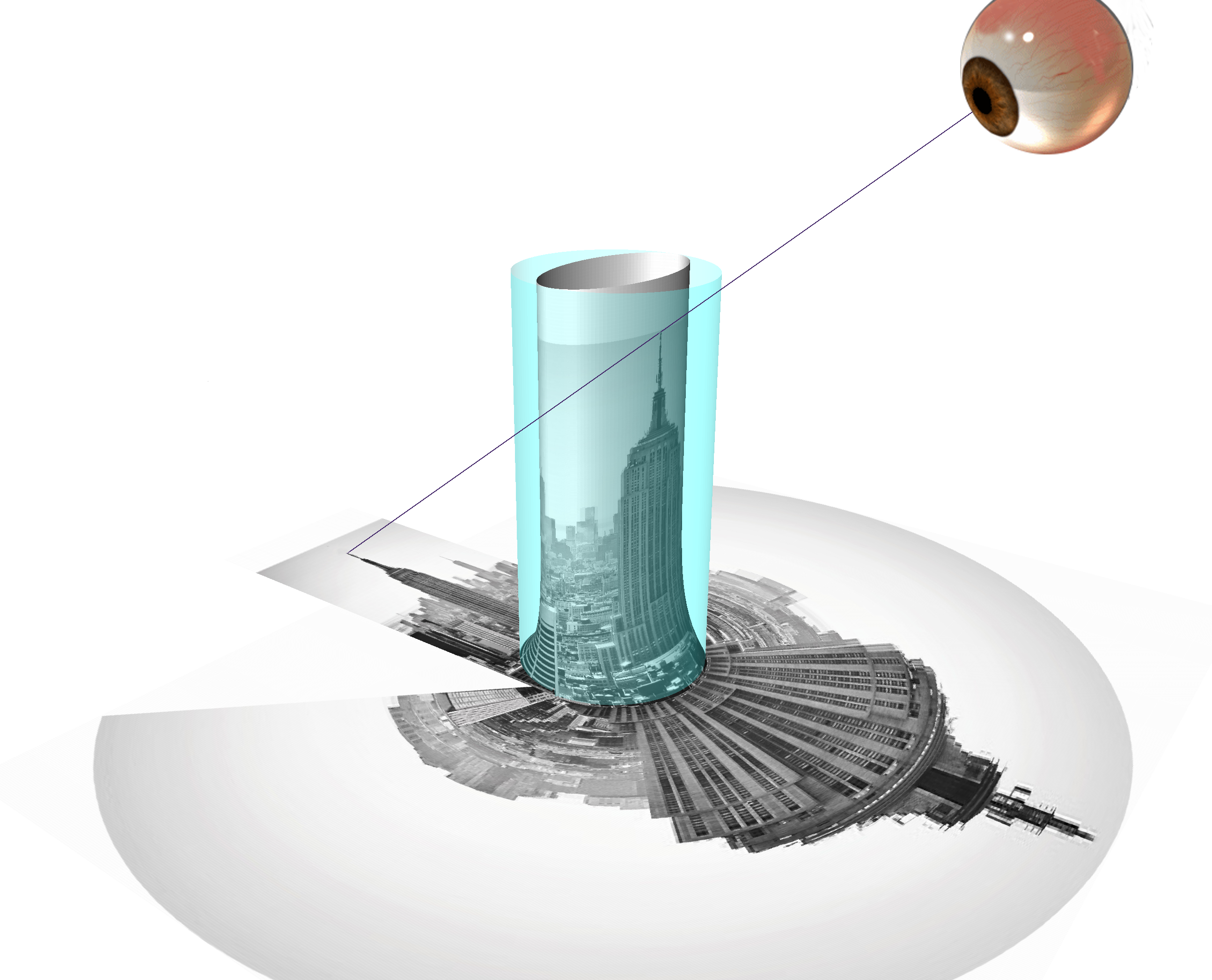}}
  \caption{The two versions of \fref{fig:correct} are visible with
    exactly the \emph{same} direction of the gaze.}\label{fig:both}
\end{figure}
A further, important, observation shows that the two
views appear along the \emph{same} line of sight, as illustrated in
\fref{fig:both}.  This implies that the two images reach
the eye as a
superposition.

The intriguing question is then whether one
can discriminate between the two superposed images.

\begin{figure}[t!]
  \CC{\includegraphics[width=0.45\columnwidth]{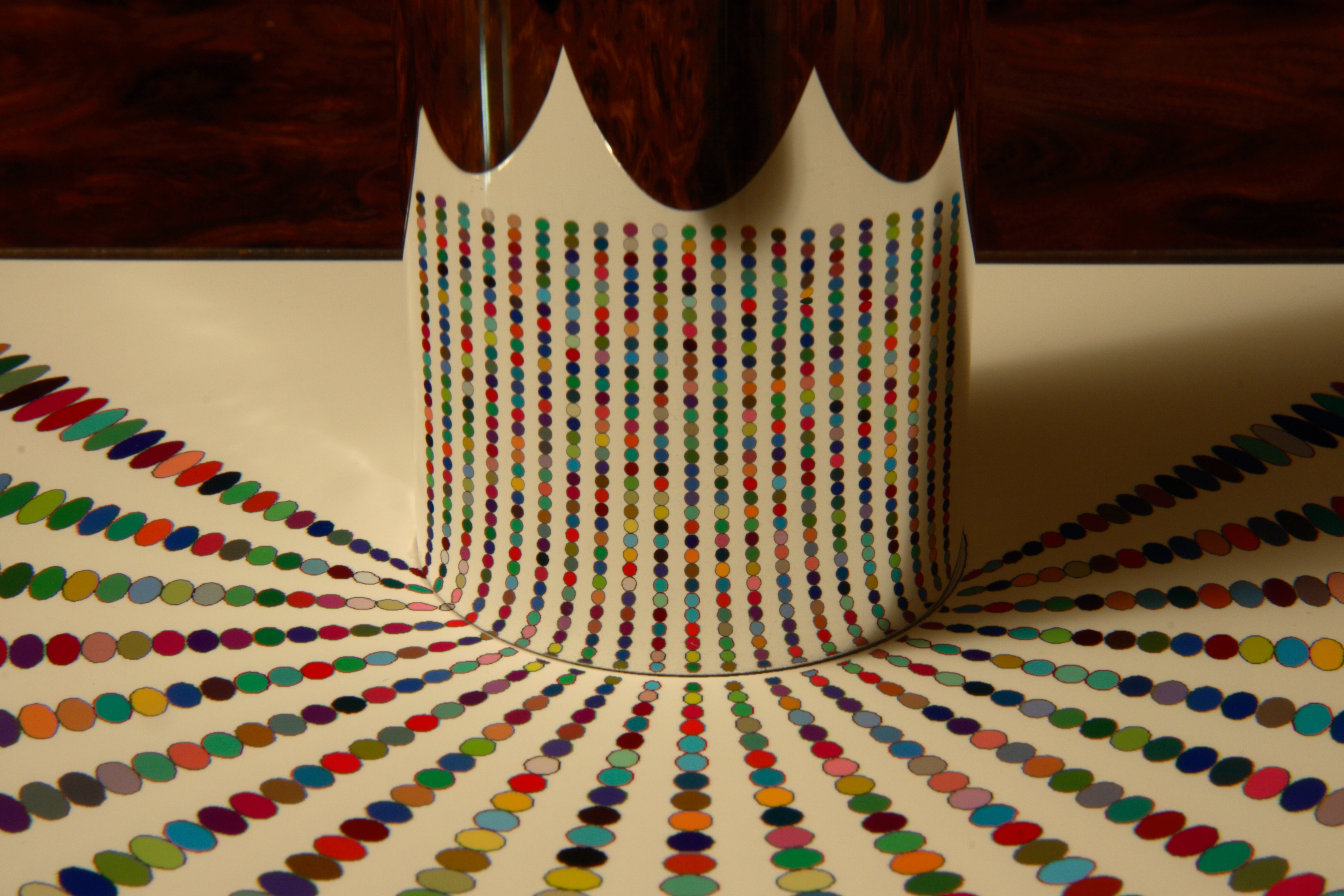}\,\,
    \includegraphics[width=0.45\columnwidth]{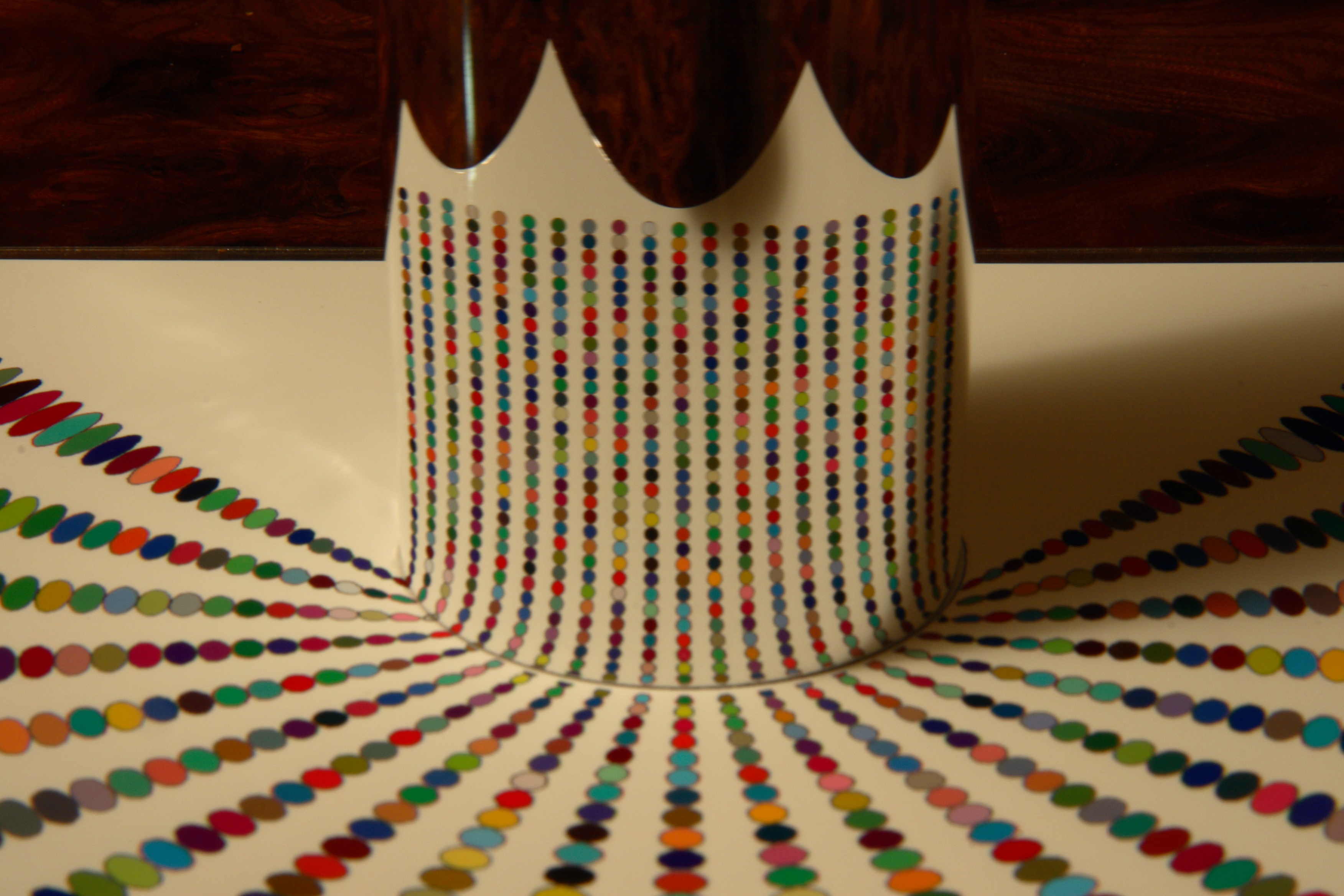}}
  \CC{\includegraphics[width=0.45\columnwidth]{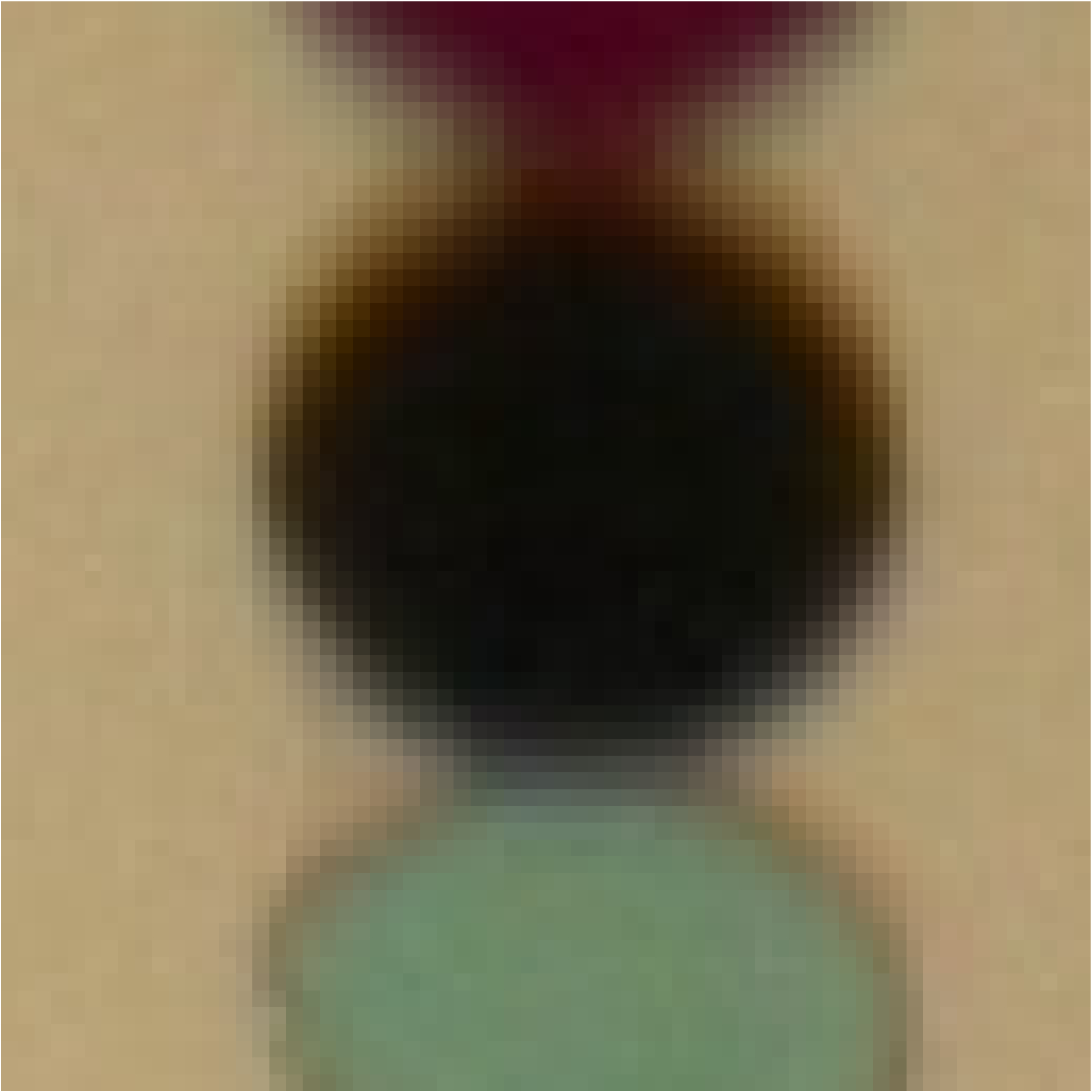}\,\,
    \includegraphics[width=0.45\columnwidth]{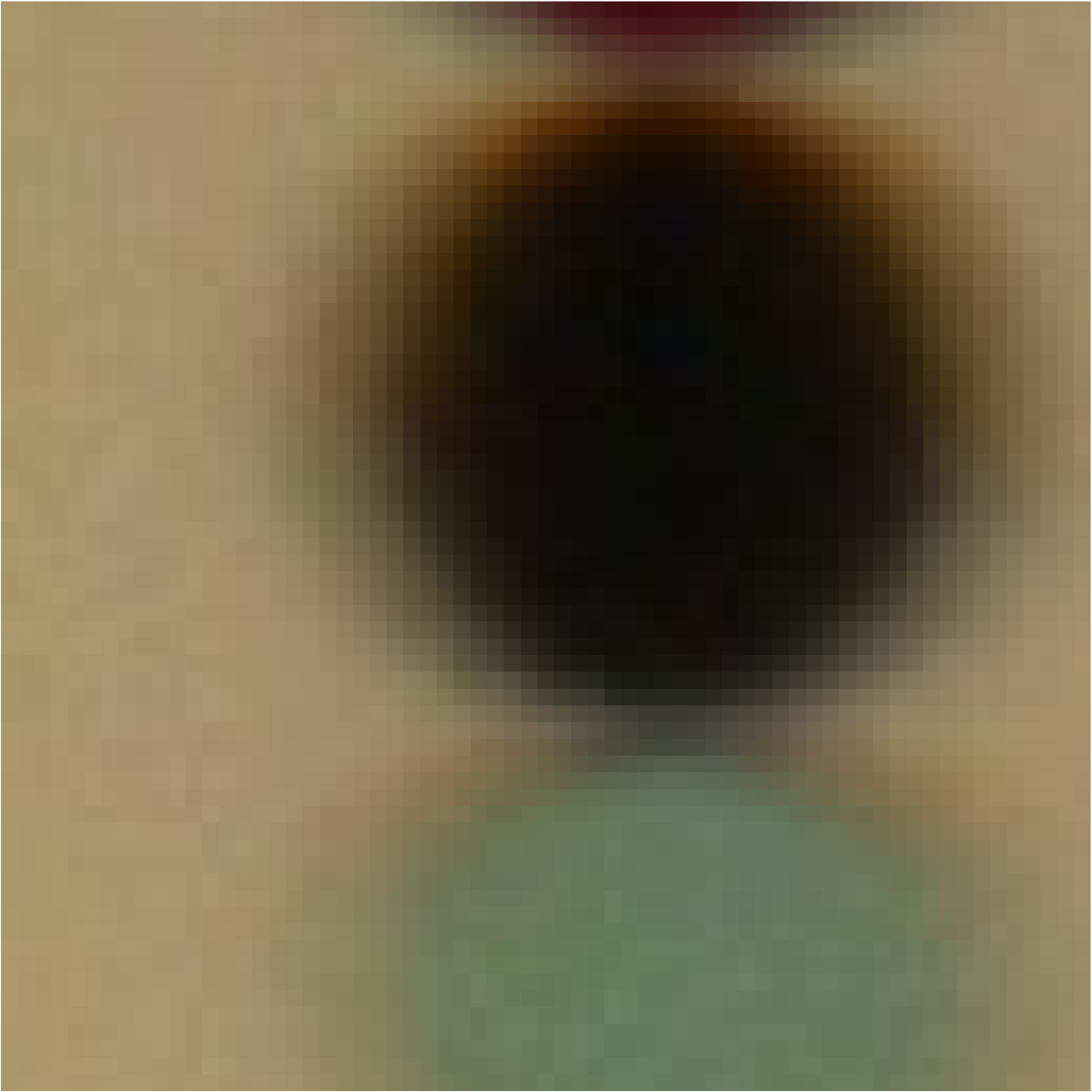}}
    
  \caption{(Photographs by Feigenbaum, August 2006) In both views the position of
    the camera and its direction of view are the same. On top: two photographs of exactly the same scene,
    reflected in the cylindrical mirror. On the bottom, magnification of a
    few dots. The only difference in making
  the two pictures is the setting of the focal distance. The pictures
  were made with a Canon EOS 30D camera and a Canon EF-S17-55mm f/2.8
  IS USM lens.  On the left, the focus is 1.2 cm in front of
  the tube's center plane; on the right the focus is  on the fifth
  horizontal
  line from the top
  of the flat image. The f-stop was f/8. The important thing to observe is the
  change in the direction of  unsharpness. On the left, the image is
  unsharp in the vertical direction, called $\HH$-astigmatism. On the
  right, the unsharpness is in the horizontal direction, called
  $\VV$-astigmatism. Note that, to the eye, the top photographs seem
  identical. (The anamorphs are of the 3D type.)
      }\label{fig:dots}

\end{figure}

We will see that this indeed is possible, and it happens in an
unexpected way. This is best understood by looking at \fref{fig:dots}.
In it, a pattern is seen on the table, which produces a regular set of
dots on the tube.
The scene is photographed with a camera, but the focal distance is
changed between the left and right takes.
Note that neither of the two choices of focal distance produces a 
sharp image, as can be seen at the bottom of \fref{fig:dots}. 
Furthermore, no other choice of focus of the camera
can make the images of the dots sharp. But the unsharpness is not
arbitrary: Both images are unsharp in a characteristic way:
One image is vertically unsharp (called $\textbf{H}$ throughout
the book) while the other is horizontally unsharp (called
$\textbf{V}$).
The letter \textbf{H} indicates that the line between the two eyes is
horizontal. Upon turning the head sideways (as explained later), the
line between the eyes will be vertical, thus \textbf{V} is used.

So the camera produces \emph{no} good image of what is perceived in
the mirror. However, the human observer perceives an image which seems
sharp, and in fact, there are \emph{two} possible sharp images to be
seen, as sketched in \fref{fig:correct}.

Since the viewer has two choices of seeing the reflection of the dots
in the mirror, the question which Feigenbaum asked is:
Which choice is preferred? It turns
out that the vertical unsharpness---\emph{vertical relative to the
natural orientation of the head}---is preferred by the eye-brain
system. That is, 
our visual system prefers the $\textbf{H}$ over the
$\textbf{V}$. This seems well-known in ophthalmology: If a patient has
vertical astigmatism (``axis'' in ophthalmological prescriptions)---called WTR (with the rule)---there
is much less need for correction than if the astigmatism is
horizontal--called ATR (against the rule). So the preference seems
somehow universal.

In \fref{fig:dots}, which is a photograph of the
image in a metallized tube, this unsharpness is clearly visible.
The uninitiated reader will not notice any
difference in the two top figures; but the eye \emph{does}. To be more
precise, any photographic image can not really distinguish between the
two possible views, as it will always record a superposition. Only the
artifact of focusing at a specific distance, as in the bottom of
\fref{fig:dots}, indicates at least some difference between the
two possible views.

Feigenbaum also shows that the
human eye cannot focus
simultaneously at the two distances. This means that the viewer must
choose (unconsciously) one of the two views, and, as I said, the
$\textbf{H}$ view is preferred. Furthermore, as the two images are in
the tube or on the table, their distance from the eye is not the
same, and this allows us to enhance the effect, by choosing where to
focus. As I said before, the image in the tube (and on the table) is
actually unsharp, but we perceive it as sharp, because the eye-brain
machinery is insensitive to vertical unsharpness (astigmatism).

The difference of focal distance is actually
more pronounced in another experiment shown in \fref{fig:pencil2}, and so
many people seem to see the effect better in that case. This one is
easy to make with a rectangular box, filled with water, and a ruler
(see \fref{fig:pencil3}).

\newpage
\onecolumngrid

\twocolumngrid

\begin{figure}[ht]
   { \includegraphics[width=0.6\columnwidth]{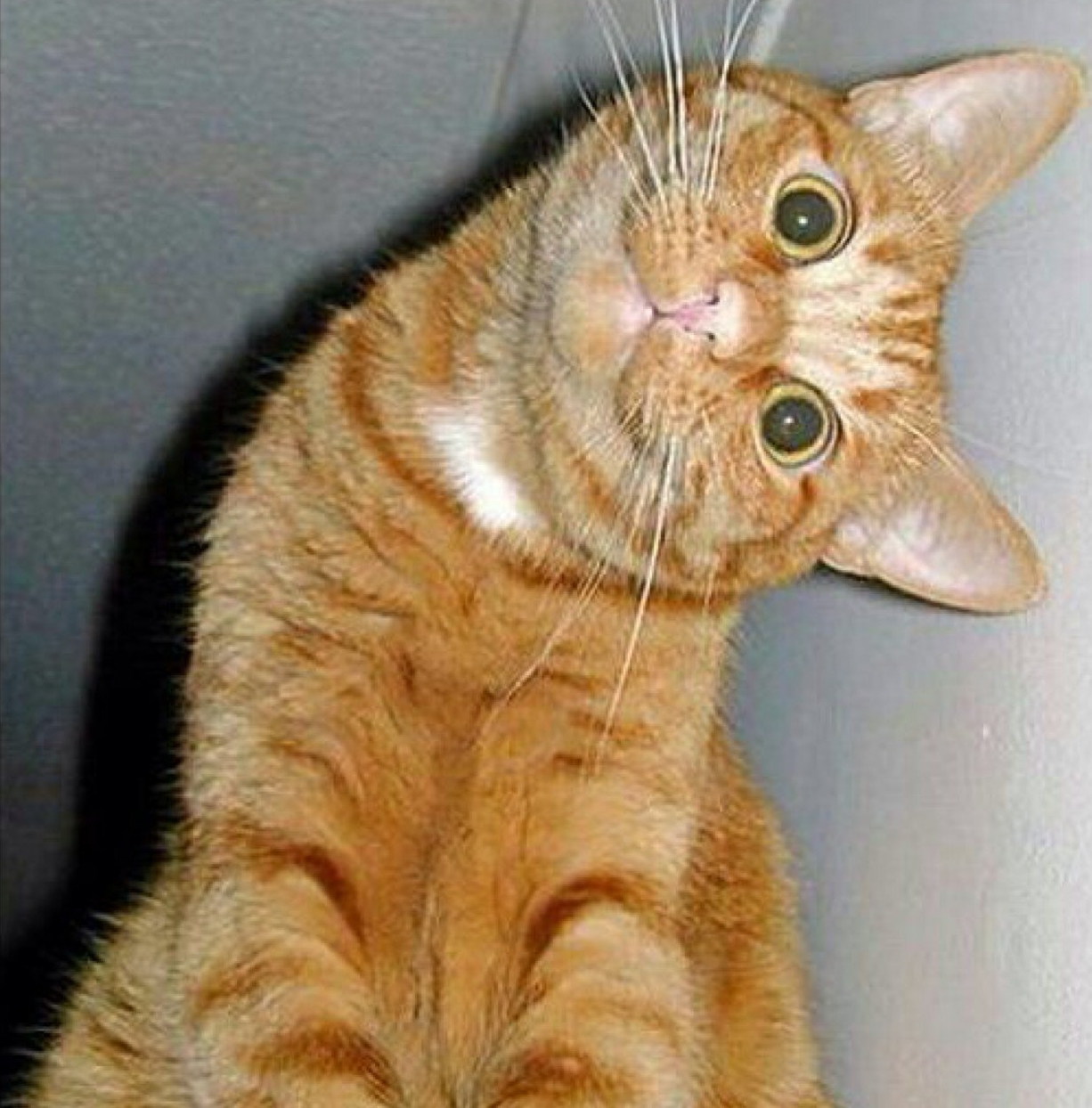}}
 \caption{Illustration of how to turn the head 90 degrees sideways. (Source:unknown.)}
  \label{fig:cat}
 \end{figure}

Since we seem to prefer the vertical unsharpness, Feigenbaum suggested
that you turn your head 90 degrees sideways as in \fref{fig:cat} (see
later for how exactly 
you are supposed to do this). Then, clearly, the notions of
vertical and horizontal get exchanged. And now, suddenly, the
\emph{other image}, the $\textbf{V}$, is going to be preferred.
You will see the reflection in a different location. With the head in
the upright position, the image appears in the tube; with the head
turned 90 degrees (and keeping the eye more or less in the same
position) the image seems to appear on the table, behind the cylinder
(it lies down). Note that the direction in which the image is seen is
unchanged, but the distance where it seems to appear depends on the
orientation of the head.

To summarize: There are two images, neither of them sharp. And there
is \emph{no} sharp image available. In such a case, the eye-brain
system will prefer the image which is unsharp vertically, relative to
the orientation of the head.

Since no image is sharp, Feigenbaum and I devised the
notion of \emph{\xx{non-object}} for what is presented to the eye.
In contrast to what is seen in a flat mirror, the non-object is not
really localized. In normal perspective, objects just present to the
two eyes two
different views of something which is fixed in space. The non-objects
are not fixed in space, so that their image seems to be located in
different points in space, depending on the vantage point. If you move
your head a little bit right-left, the image seems to turn around the cylinder.

(This unsharpness of the non-object has the amusing consequence
that \xx{autofocus}ing with digital cameras will be confused by the
two unsharp images which are at different distances, as there is no
ideal focal distance available.)

What is the reason for the possibility of seeing two choices?
As mentioned before, the eye is not a pin-hole camera, but has
some \xx{aperture}. This simple, but important fact means that
different points on the \xx{retina} see different caustic points
(points of maximal intensity, see ``The main ingredients,'' below). In
other words, to really understand what one sees, one has to consider
the image, the tube, \emph{and} the eye. This is what Feigenbaum
did.

\ssection{The main ingredients}\label{sec:mainingredients}

The study of Feigenbaum combines in a clever and certainly completely
novel way certain observations which are well-known to people
familiar with optics. I will explain now what the pieces are, and how
they are combined.

\begin{figure}
  \CC{\includegraphics[width=\columnwidth]{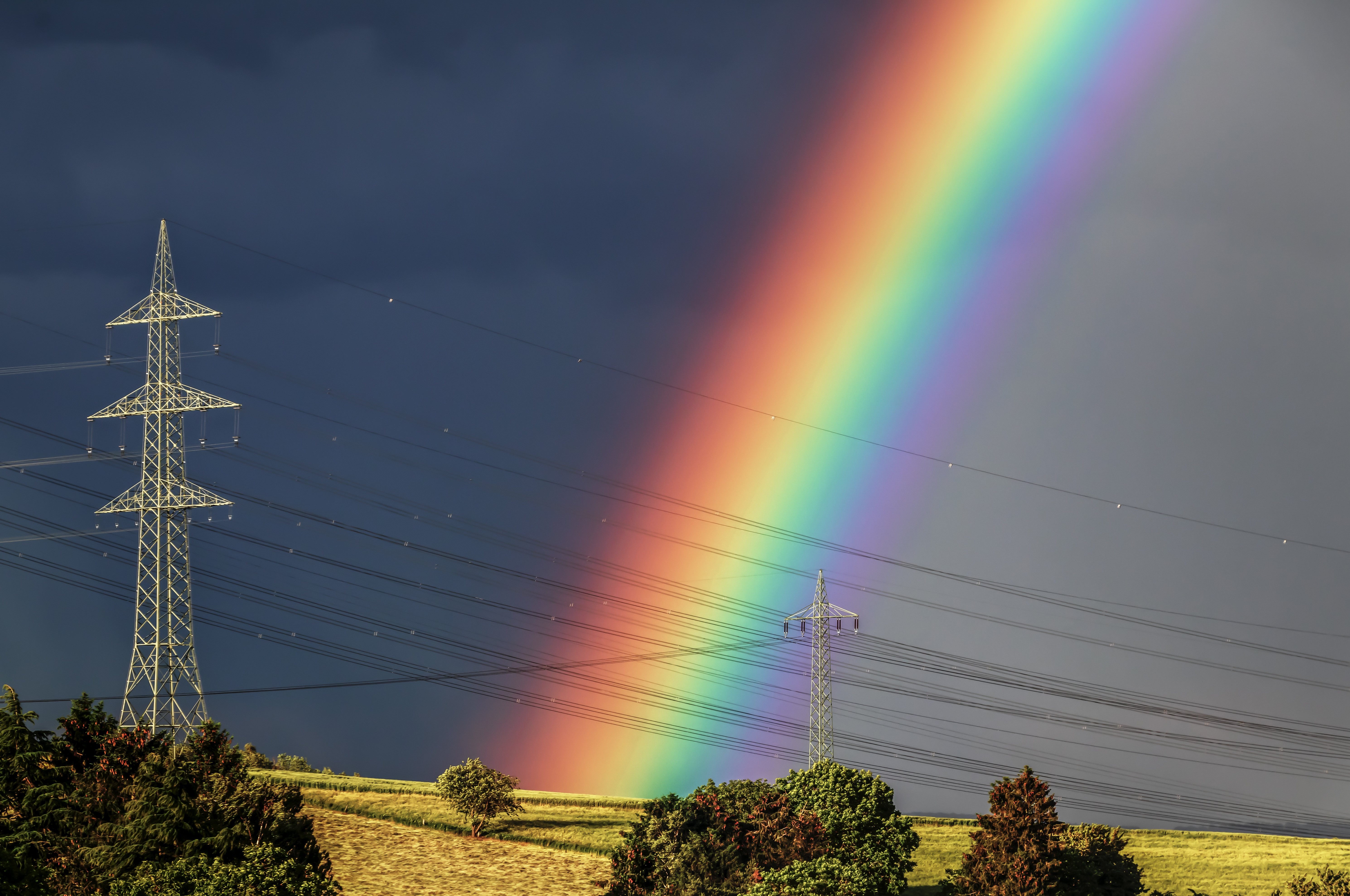}}
  \caption{A typical rainbow. Note that the sky on the outside of the
    rainbow is darker than it is on the inside of the rainbow. (Source:https://www.maxpixels.net/Rainbow-Meadow-Sky-Strommast-Landscape-Nature-4285843.)}\label{fig:genevarainbow}
\end{figure}

\begin{figure}[ht!]
  \CC{\includegraphics[width=\columnwidth]{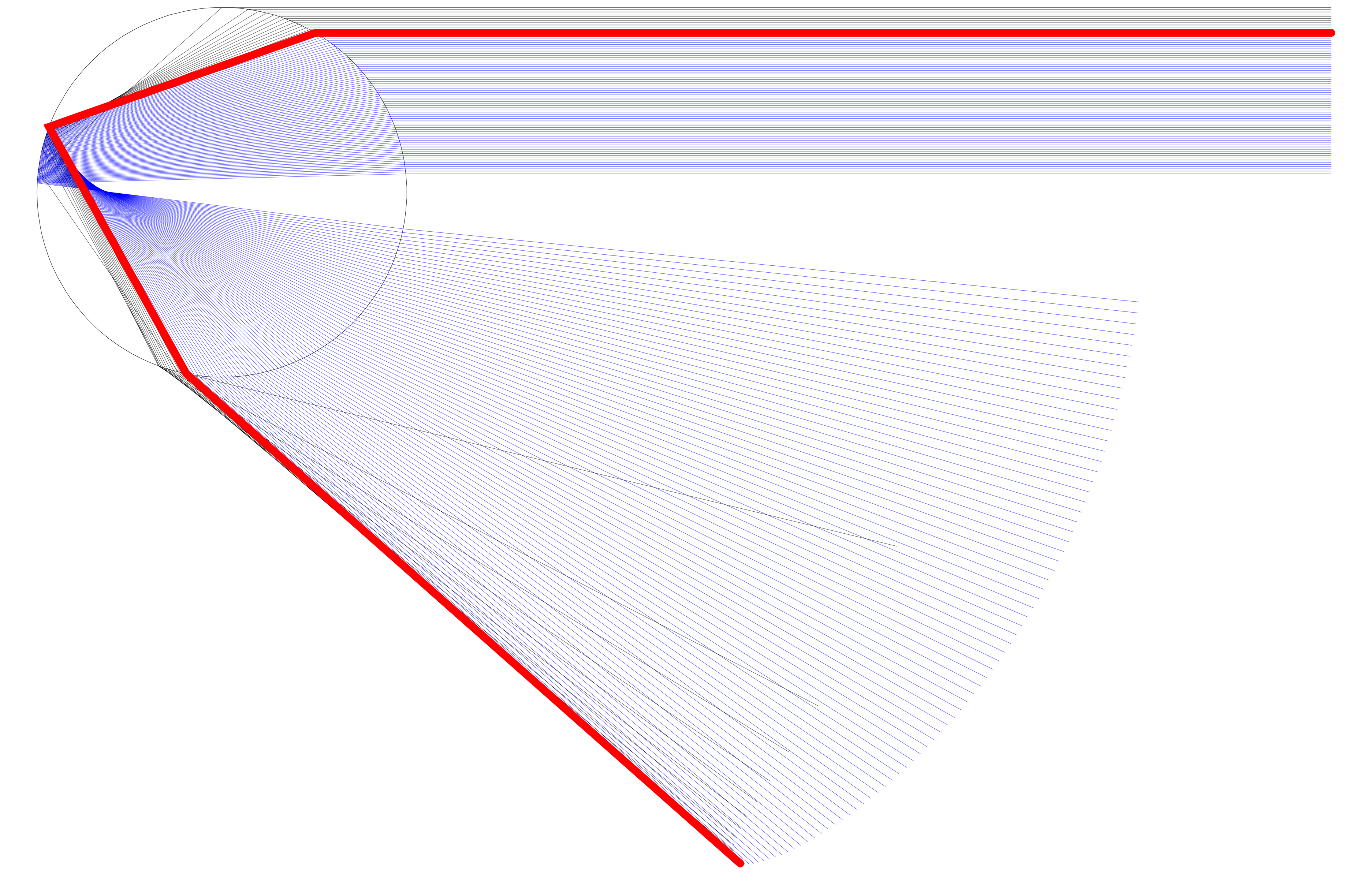}}
  \caption{Illustration of the caustics for a rainbow (for just  one
    color/wavelength) in a spherical
    (actually cylindrical) raindrop. The rays from the sun reach the
    drop completely parallel from the top right. They get reflected
    once on the interior left side of the drop, and leave to the right. One can
    see that many incoming rays get dispersed, but many concentrate
    near the red line, from both sides, blue and black rays. This
    gives more intensity, along the red direction. The
    envelope to which the red line is tangent on the lower right is
    called the \emph{caustic}. This sketch also explains why the
    exterior of the rainbow is lighter on the inside: This is because
    the dispersing rays come out at a flatter angle (and will be seen
    in drops which are ``below'' the rainbow). For more complicated reflection
    patterns,
    I like to look at the beautiful photographs  in the classical book
    \citep{minnaert}. For an accessible theoretical discussion, see
    \citep{nussenzweig1977}. The Moir\'e artifacts in the image come
    from the pixeled nature of reproducing images.}\label{fig:rainbow}

\end{figure}

The first ingredient is the role \emph{\xx{caustic}s} play in vision.
Every lay person has probably seen caustics without realizing it,
because they are what one
sees in rainbows. What people generally do not know is that the
small drops, which form the lenses for the rainbow, are quite terrible
lenses. Indeed, as shown in \fref{fig:rainbow}, the light rays which
enter the drop do not come out of the drop at a sharp angle,
but are rather fanned out. This raises the question of why the rainbow
seems reasonably sharp when we look at it. What happens is that the
density of the rays in the fan is not uniform and there is one
direction at which the density is highest. And this is where you see
the rainbow, each color at a different angle, but quite sharp. The
conclusion to draw is that we do not see the spread-out rays, but
rather these
highest density regions which are the caustics.

Another observation you should make is also
seen in  \fref{fig:genevarainbow}: namely, the sky around a  rainbow seems
to be dark on the outside and 
diffuse on the inside, and this is explained again by
\fref{fig:rainbow}, which shows that the fan opens only upward, but
not downward.\footnote{The reader might be confused by this as in
  \fref{fig:rainbow} there are only rays fanning upward from the
  caustic. But a little thought shows that we actually look at this
  figure from below and see images of reflections in drops at
  different heights.}
The conclusion to draw from this discussion is that the correct way to
study what one actually sees must go through a study of caustics, and
where they are. In particular, the usual 
\xx{ray-tracing} methods often found in 2d-projection graphic programs
are not adequate. In the case of reflection from a tube, ray-tracing
would do the following, as illustrated in \fref{fig:erect}. Take any image, and wrap it onto the tube. Fix
the position of the observer, and assume the eye is just a point. Draw
a line from this ``eye'' to an image point on the cylinder, and
continue it down to the table, using the rules of reflection from the
cylindrical mirror. This is the ray-tracing image of the scene on the
mirror. It is \emph{not} the way Feigenbaum constructs his anamorphs,
because, as shown in the left side of \fref{fig:correct} one of the
two images appears on a surface whose cross-section is an ellipsoid (of
ratio $2:1$ at infinite height), and so one must wrap the image somehow on this surface
(and not on the surface of the tube).

The second ingredient which is one of his main insights
is the question of what happens if
there are \emph{two} caustics presented to the eye.
As I alluded to on page~\pageref{pag:2d}, there are, most often, two
caustics to be seen because any mirror is, in every point, curved in
two principal directions.\footnote{This is not to be confused with
  what earlier authors call two caustics, namely the two curved pieces
  of the 2-dimensional \fref{fig:more1}
which make what here is called one caustic only.} (If the mirror is completely flat, the two
caustics will coincide.) The tube is a particularly nice example
because it is vertically flat, and horizontally just a circle, which
allows for an explicit calculation of the caustics. In the paper
\citep{witchball}, the interested reader can find a variant of
Feigenbaum's calculation, where the two caustics are determined for
reflections from a sphere.

While the book explains
that the phenomenon of two caustics is ubiquitous, the most striking example is
``reflection on a tube,'' the title of this book. The setup is shown in
 Figs.~\ref{fig:anamorph} and \ref{fig:istvan},
 and in many illustrations throughout the book.
 \medskip
 
 \emph{After discovering for the first time that there are two
   disconnected caustic
   images for the cylinder, the interesting---and in my view, completely novel---question
   which Feigenbaum asks (and answers) is
which of the two options is preferentially chosen. Furthermore, an
explanation is given for how  this choice is made. And, as will
be seen, this finally must be related to how the visual system of
animals processes the inputs it gets. }

 \medskip
 
There are several precursors in the literature on viewing objects in
water, where the authors knew that there are two images, with different
astigmatic directions. As far as I can tell, while some mention a
binocular effect, it seems that the cyclopic effect has not been
mentioned. See, \eg \citep{Kinsler1945,Bartlett1984,Horvath2003}.
(For a review on classical caustics, see \citep{Berry1980} and the
original work on ``catastrophe theory'', \citep{Thom1972,Thom1976,Arnold1974,Arnold1975}.)


\ssection{How to look at the cylinder, that is, the tube}

It is easy to make your own cylinder. Best results are obtained for
cylinders of diameter about 5 cm (or 2 inches).
The cardboard kernel of a kitchen-tissue roll is just about
right. One wraps a reflecting (silvery) Mylar sheet tightly around it
and places the whole thing in the center of the included figure. (Do
not use aluminum 
foil, it is not flat enough.)

To make an anamorph of a jpeg file of yours,  go to the page
https://fiteoweb.unige.ch/$\sim$eckmannj/a4\_shift.html\label{pag:html}\ which also contains the necessary instructions of how
to print the new jpeg which it constructs. It is important that the
size of the printout is correct.\footnote{I thank No\'e Cuneo for
  transforming Feigenbaum's program (which was written in pascal) to
  the html version.} Recall that there are 3 possible anamorphs, called
``erect,''
``flat,'' and
``3D''.
The program will produce the 3D version, which
distorts the image minimally when seen inside the tube. This is the
one which is most natural, and is seen in the \textbf{H} direction. The flat one
is undistorted when viewed in the \textbf{V} direction, and the erect one is
the ray-tracing. (Note that there is no choice of anamorph which is
undistorted in 3D and flat simultaneously.)

The table should be
as flat as possible. The eye of the observer is supposed to be
positioned at a distance of 
about 25 cm from the cylinder and about 40 cm above the  table. In this
position, most people see the image of the drawing as if it were
inside the cylinder. Closer inspection, and the calculations by
Feigenbaum, show that the image appears not on the surface of the
cylinder, but on an ellipsoid half as thick as the cylinder.
Now one should turn the head 90 degrees sideways, 
but in such a way that the eye with which one looks (you should look
with 1 eye, see below) remains in the same place.
In this configuration, most people see the drawing ``lying
down,'' as if it were reflected behind the cylinder. (This rule of how
to turn the head is important because you should, as in
\fref{fig:both}, not change your line of sight.)

It is important to
note that the effect has \emph{nothing to do with \xx{binocular vision}}, as you
can check by covering one eye. However, as Feigenbaum studies in
detail, binocular 
vision into the tube is quite different from binocular vision of true
objects, because the two eyes see two different non-objects. (This is
then an over-determined problem for the brain.)

I showed, in 2009, Feigenbaum's project
to a neuro-ophthalmologist friend of mine,
 \xx{Avinoam Safran}. He pointed
out that turning the head is less good than turning the cylinder (and
the ``table''). In fact, the inner ear signals the position of the head
to the brain, and he told me that the 4th cranial nerve activates a
muscle which rotates the eyes towards the nose.\label{pag:psycho} ~I encourage any
interested reader to do the experiment in this more
complicated, but cleaner way. Another possibility would be to do this
in the space station, where gravitational orientation is
missing. Finally, I decided to have a hologram made, which 
avoids these problems,  see page \pageref{pag:holo}.

The book also explains where exactly the image appears. If you
watch closely, (with both eyes), you will notice, as I said before,
that the reflected image
appears glued onto an {elliptic surface}, inside the
cylinder. Furthermore,
this surface
enlarges towards the bottom to reach the circle where
the cylinder touches the table. (The reader can see this illustrated
on the left side of \fref{fig:correct}.
The alternate image, however,
is completely flat on the table (right side of \fref{fig:correct}).

Another experiment, related to the different astigmatisms of the two
caustics, is to move your head slightly up and 
down, the picture inside the tube will move with this vertical motion, sideways motion of  the head
makes it turn around the cylinder. The roles are exchanged if you focus on the image
``behind'' the cylinder.

\ssection{Other examples of multiple caustics}

The aim of Feigenbaum's work is to shed light on the questions these
observations raise. I just list some of the issues
before I explain some of Feigenbaum's further contributions in his manuscript: 
\begin{itemize}
\item Is the ``pixel'' \xx{resolution} of the eye (\ie its \xx{acuity}) good enough to
  actually distinguish  the 
difference  between the two non-objects, \ie between the 2 images of \fref{fig:dots}?
\item When one uses both eyes, the two eyes get two different images
  from the same scene. What is the magnitude of the effect that
both eyes get slightly different images, and how does the brain
interpret what is seen if the image is coming from a non-object?
\end{itemize}

Such questions lead to  a discussion of the quality of human eyes versus
fish-eyes, which are shown to be vastly better than the human
eyes. From this Feigenbaum draws some conclusions about  the
relatively bad evolution from fish-eyes to the eyes of
land-animals. 

\begin{figure}
  \CC{\includegraphics[width=\columnwidth]{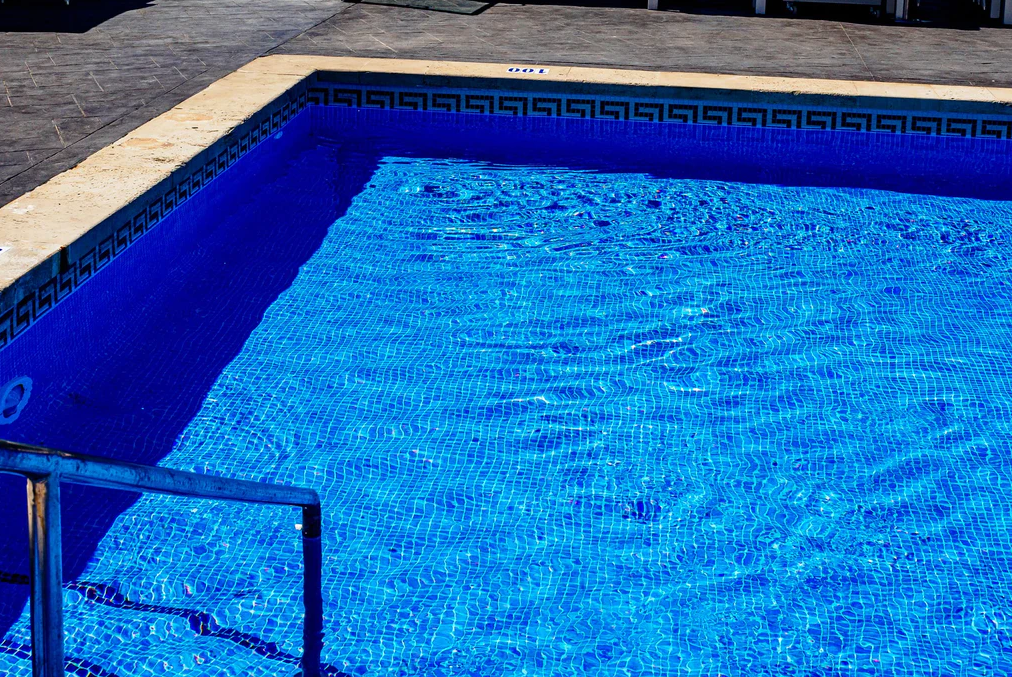}}
  \caption{The up-sloping pool floor. In this photograph, the pool
    floor clearly seems to go up towards the far corner. I could not find a
    photograph which also shows that the slope gets flatter as you
    look far away, as sketched in \fref{aafig5}. (Photograph: From the
    online catalog of Oasistile.)}\label{fig:pool}
  
\end{figure}

\begin{figure}
\CC{\includegraphics[width=\columnwidth]{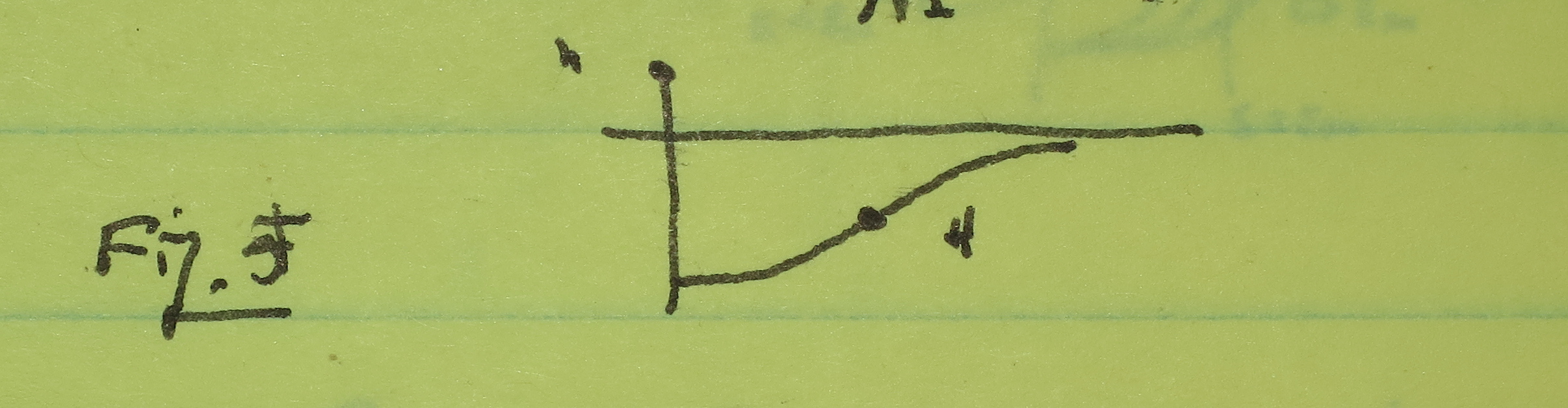}}

    {\caption{The apparent bending of the sea floor for an observer
  standing at height $h$ (on the left) with the caustic point $H$
  indicated. If the viewer looks down at an angle of $35$ degrees, and
the pool is $10'$ deep and the viewer's eyes $10'$ above the surface,
the floor will seem to slope upward by about $10$ degrees.
      } 
\label{aafig5}}

\end{figure}

\begin{figure}
  {\includegraphics[width=0.5\columnwidth]{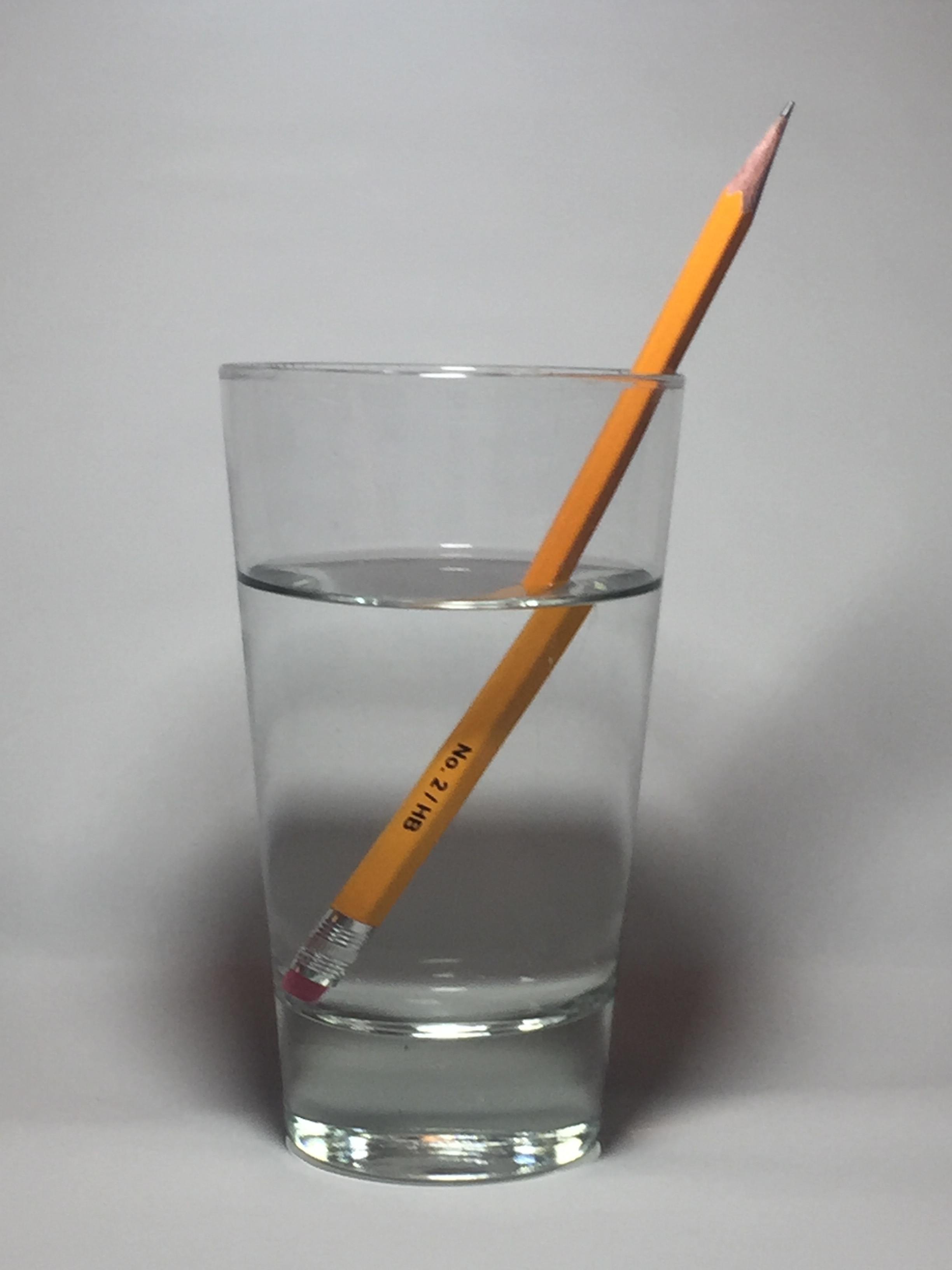}}
  \caption{A standard photograph of the ``broken'' pencil. However, as
    shown in \fref{fig:pencil2}, Feigenbaum was interested in what one
    sees when looking into the water \emph{from the air}, at a shallow
    angle. This standard image is a view \emph{through} the water.
   }\label{fig:pencil} 
  \end{figure}

\begin{figure}\label{pag:pencil}
  \CC{\includegraphics[width=\columnwidth]{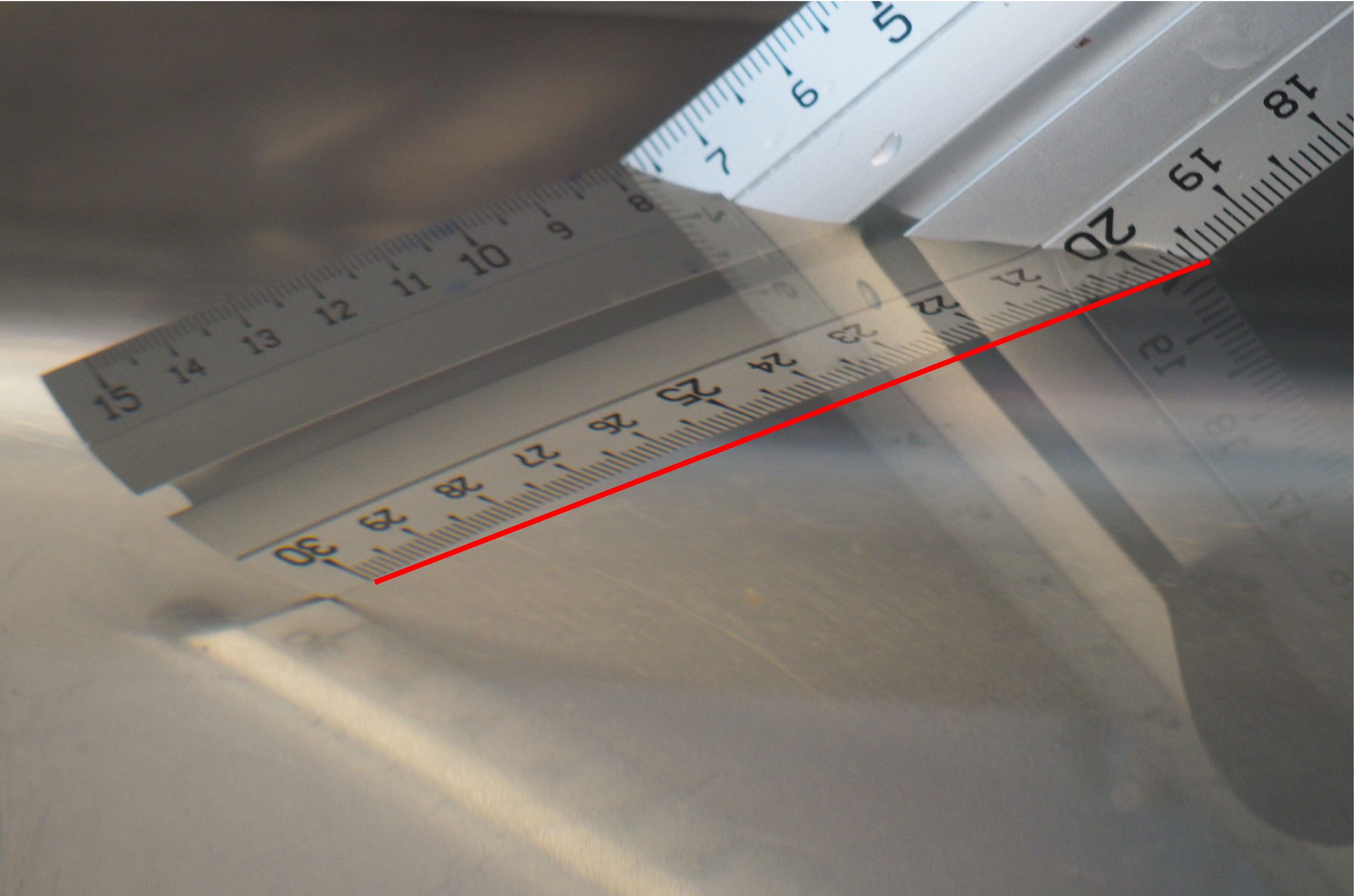}}
  \caption{A photograph of the ``broken'' pencil (actually a ruler). The point Feigenbaum
    made,
    is that the part of the  pencil \emph{in} the water,
    when looked at through the top
    surface of the water, at a very low angle, is \emph{not} straight. This
    is especially well visible
    at the lower edge of the ruler, and was sketched by Feigenbaum.
    In \citep{ruler}, I worked out the details of what Feigenbaum had
    in mind. 
    Note that the ruler is strictly perpendicular to
    the line of sight from the camera. If you repeat this experiment
    and tilt your head 90 degrees, you will see that the left bottom
    corner moves toward you. The left top corner will also seem to
    move toward you, but less, and so the whole ruler seems to
    rotate toward you.  The interested reader can also see the effect
    of the ruler in a setup as in \fref{fig:pool}: In that case, the
    vertical bar of the handrail which goes into the water seems bent
    towards the viewer, and the distance of the bottom end will change
    if the viewer rotates the head by 90 degrees. (The effect is
    strongest if the eye is close to the surface of the water.)}\label{fig:pencil2} 
  \end{figure}

\begin{figure}\label{pag:pencil3}
   {\includegraphics[width=\columnwidth]{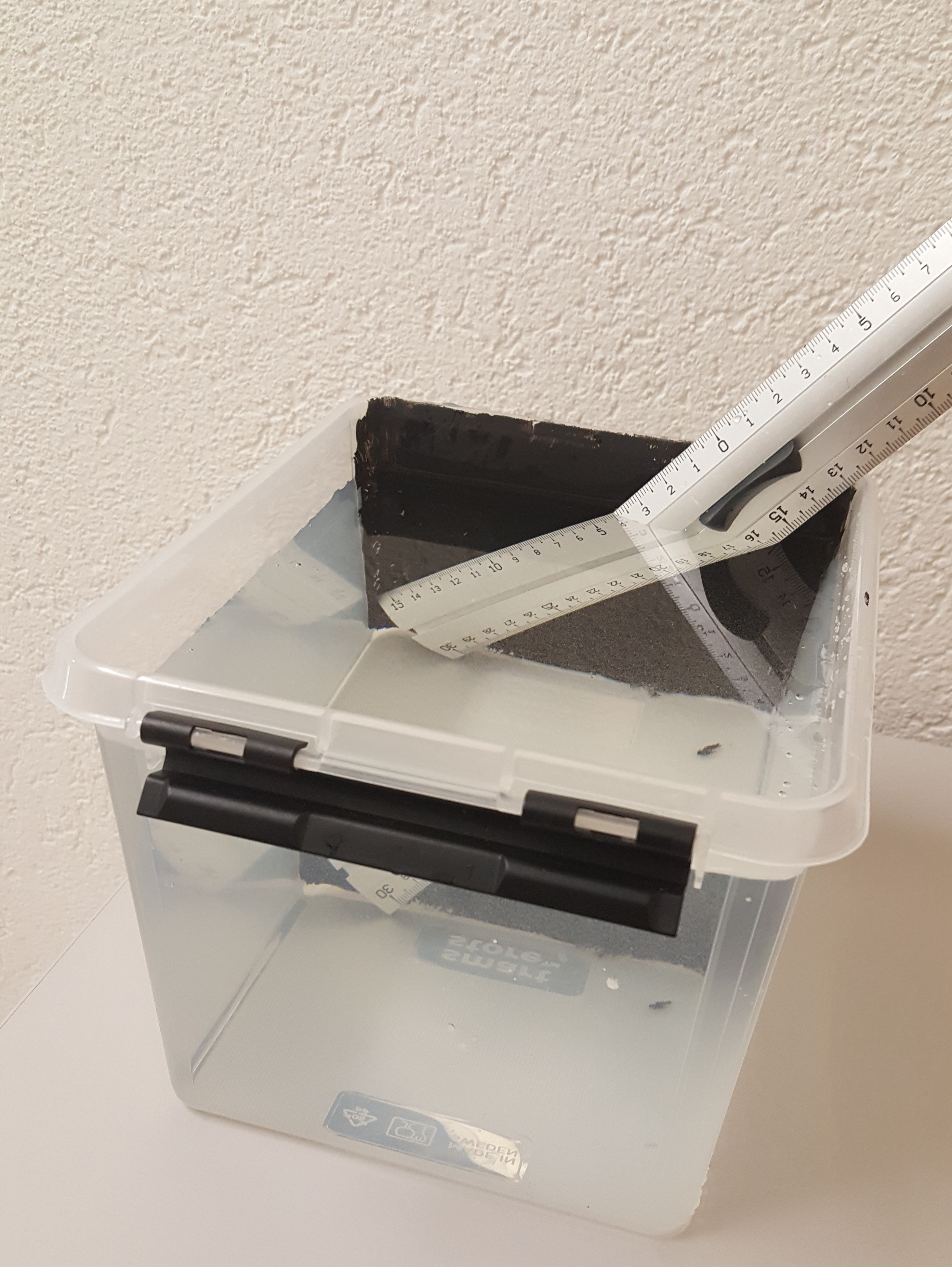}}
   \caption{A convenient setup for the ruler. It is best to fill water
     to the rim and to look into the water from above at a flat angle.}
   \label{fig:pencil3}
  \end{figure}

\begin{figure}\label{pag:archer}
     \CC{\includegraphics[width=0.8\columnwidth]{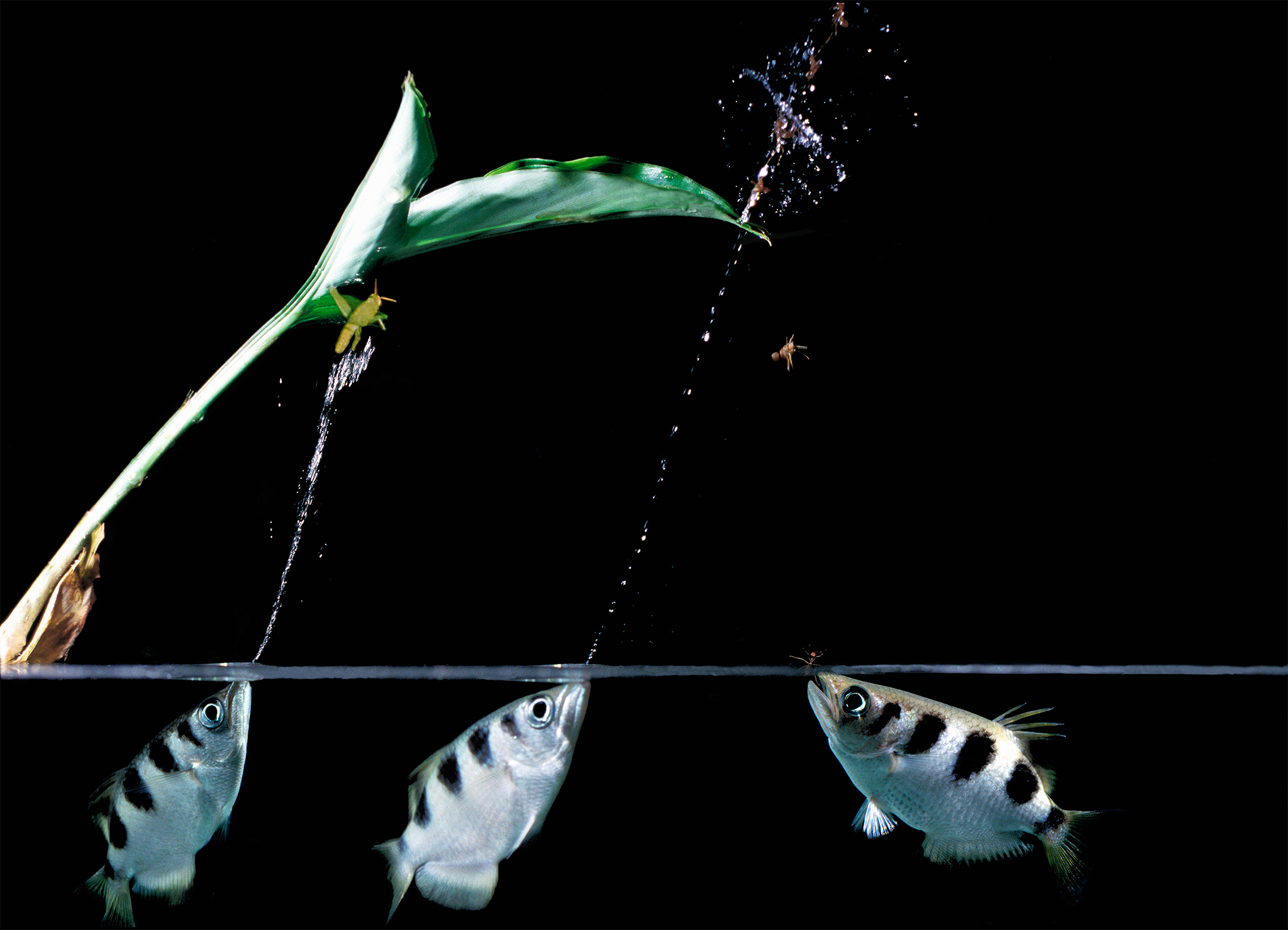}}
     \caption{Illustration of what the archer fish can do. Note that
       the eye is below the water, and therefore, the fish must
       ``calculate'' not only---because of the different indexes of
       refraction of water and air---in which direction to aim, but also where the object will
       fall. Feigenbaum has a section in which he does the
       calculation for this case, which is a variant of the cylinder
       case. He devotes a chapter ``Caustical Imaging at the Air-Water
       Interface'' on this.
       Archerfish of the
       species Toxotes jaculatrix take down insects in
       Indonesia. Photograph by A\&J Visage, Alamy, National
       Geographic, September 4, 2014.}\label{fig:archer}
       
   \end{figure}

The appendices contain furthermore several beautiful examples based on
the \textbf{H}-\textbf{V} dichotomy where caustics
appear in everyday 
life. I illustrate these in \fref{fig:pool}--\ref{fig:archer}, and
some explanations are summarized in the captions:
\begin{itemize}
\item Looking from air into the water, \emph{from above} (a pool, the
  sea) as in \fref{fig:pool}.
\item The bent (broken) pencil as in \fref{fig:pencil}--\ref{fig:pencil3} (one can see again \emph{two}
  views).\footnote{One of the two views was certainly known at the end of the
    $19^{\rm th}$ century. I am not aware of any discussion of the
    second possible view. An example can be found in
    \citep[Fig.~316]{watson1907}. I worked out some more details in
    \citep{ruler}.}
\item Looking from water into the air: this is the problem of the
  archer fish which ``shoots'' at targets in the air from below the
  surface of the water: see \fref{fig:archer}.
\end{itemize}
It is perhaps useful to add here one of many Feigenbaum sketches for
this upsloping effect in \fref{aafig5}
with my own caption added.


\begin{figure}[h!]
   {\includegraphics[width=0.7\columnwidth]{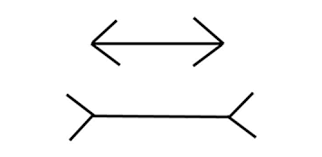}}
  \caption{The well-known optical illusion which makes the lower line
    look longer than the upper one. Apparently invented by Francis
    Xavier Hermann M\"uller, 1889. }\label{fig:arrow}
 \end{figure}

\begin{figure}[h!]
  \CC{\includegraphics[width=\columnwidth]{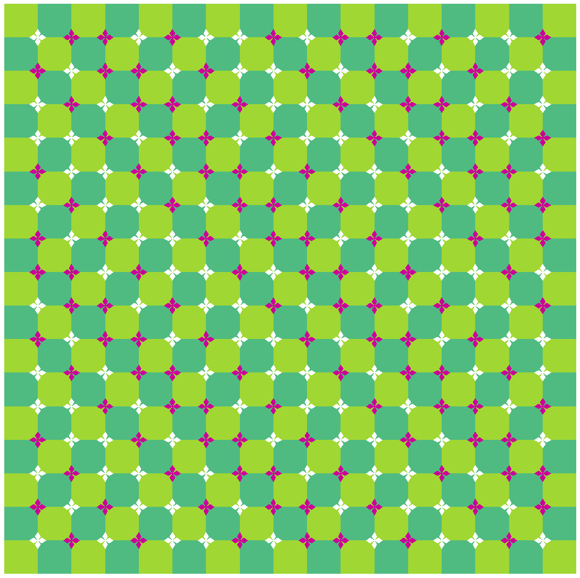}}
  \caption{A moving optical illusion: The pattern seems to move if the
    head is gently moved up and down. I guess this comes from changing
    focus of the eyes. Copyright: Akiyoshi Kitaoka, Ritsumeikan
University.}\label{fig:akiyoshi}
  
\end{figure}

\begin{figure}
 {\includegraphics[width=0.7\columnwidth]{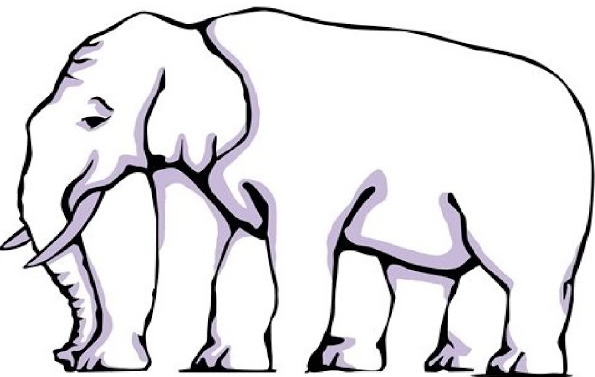}}
  \caption{The famous elephant riddle, which confuses about the number
    of legs. (This is one of many adaptations of \citep{rogershepard1990}.)}\label{fig:elephant}
  
\end{figure}

\ssection{An important comment} As follows from the previous sections,
the phenomenon Feigenbaum described is \emph{not} about
optical illusions. To see the difference, let me just show, in
\fref{fig:arrow},
a favorite
illusion
everybody has seen. Here, the effect is that the middle
segment looks longer in one case than in the other, and this is
triggered, as in many other examples of this kind, by what surrounds
the central line. The illusion disappears if one draws  vertical lines
connecting the tips of the arrows.
Other optical illusions, such as \fref{fig:akiyoshi}, exploit
where the eye puts its focus. The apparent motion of the picture
disappears if you fix your eyes at a fixed center of the picture. Finally, those like
\fref{fig:elephant} simply present confusing realities, and the puzzle
disappears if one analyzes the way the legs are drawn.


The reflections in the
tube are of a completely different nature. Because the mirror is
2-dimensional,
the viewer is presented with two
options (the vertical and the horizontal picture), and it is somehow
the brain-eye system which has to make a choice: Namely your visual
apparatus prefers
to see that picture which is unsharp in the current vertical axis of
the eye, which depends on the position of your head. No amount of
rationalizing what is seen can make the dichotomy disappear, because
there is simply no adequate best middle ground between the two pictures.

\ssection{The origin of the book}

Mitchell Feigenbaum passed away on June 30, 2019.
His interest in anamorphs was originally raised in discussions with Kenneth
\xx{Brecher}, in March 2006. By the time of his death,
he had worked, continuously, for
about 13 years on this project. Many versions of TeX files
for the book were created over this period.
Several of his colleagues  discussed with him his project,
helped typing and gave
input. Perhaps I happened to be the most insistent and enthusiastic
follower  of his project. The state of the book at Feigenbaum's death
was close to finished, but, of course, some things were missing. The
last few appendices were incomplete. Feigenbaum 
certainly wanted to rewrite the preface, which, in the form I had the
manuscript\footnote{``Manuscript'' means here the collection of TeX
  files and figure files.}, did not
refer to these last sections. While he had clear ideas of what needed
to be done, his health problems did not allow him to finish the task.
%
Given Mitchell's
huge investment of time and energy, and the originality of his findings, many of his friends felt it was
a major loss to the community if the book project were not completed, as best as possible.
However, for the reasons mentioned above regarding publication rights, this has to date not been possible.
The present article is an attempt to provide the broader community with some understanding of what Mitchell accomplished.
While it is still preferable to produce the book in its amended form,
I recognize that Mitchell would probably not have been satisfied with
what I have done.

Let me be clear about
the following issue: I did not undertake this project as an 
historian of science,
and I do not intend to guess what other thoughts Feigenbaum might have
had.
One exception to the principle of not being a historian
is my treatment of a section ``Evolution and Design'' which had
no text. Since it is an important point, I have added some correspondence from Feigenbaum, see page~\pageref{pag:ead}.

\medskip

\emph{I think the book should be understood as to how a study of the optics
  of vision is closely related to questions of how the visual system
  and the brain interact with the information which comes through the eyes.}
\medskip

My experience is that his unpublished book, as well as other publications by
Feigenbaum \citep{feigenbaum1978,feigenbaum1979} will need some adjustments by the reader.
There are two related reasons for
this:

First, Mitchell's language is, as noted by John Horgan (in ``The End of
Science'' (Interview 1994) \citep{Horgan1997}) 
``as if English were a second language he\footnote{That is, Feigenbaum} had mastered through sheer
brilliance.''

The second aspect which makes the reading not so easy, is
Feigenbaum's technique of proving statements. While many scientists
are willing to use sentences like, ``the following calculation is left
to the reader,'' or to cite a reference to a known calculation, Mitchell would
never allow himself such a liberty, which means that he verified everything
from scratch and gave all details. But of course, some training is
required to know which details to gloss over on a first reading of the text.
James \xx{Joyce} is also
difficult to read but, in the end, that makes for a rewarding experience.

I attribute Feigenbaum's style to what I like to call a 19th century mind:
He was not only extremely critical of others, but even more
rigorous with himself. And so, as I said above, I am not sure that my friend
Mitchell would have been happy with how I would finish the book or even discuss its
highlights as in this introduction. But I think
it is important that those parts of the text which were complete at
his death are left ``as is,'' so that they reflect his personal taste
and style. As I had checked, criticized, and discussed with him all
calculations before 
his death my task was reduced to adding missing items,
where possible, and giving some explanations we discussed earlier.\footnote{There were about 50 unresolved questions of
mine, and Sects.~I and J were incomplete.}

\ssection{Feigenbaum's interest in vision}

Let me now come back to the science in Feigenbaum's book, and how it
is related to his general outlook.

Many people are familiar with  Feigenbaum's discovery about the
\xx{universality} of \xx{period doubling}, and the constant
$4.6692\dots$ which now takes his
name, \citep{feigenbaum1978,feigenbaum1979}.

Obviously, one  may think that this book is an addition to his
very successful findings of
the 1970's. But he
developed here 
a very different subject, namely optics, and in particular its
relation to vision, the construction of the eye, and the
interpretation by the brain of what is being seen.

Feigenbaum had, for many years before he started to write this book,
been interested in optics and vision. This began with work on the
problem of why the moon seems to be larger when it is near the horizon
than when it is high in the sky. This appears to be unpublished. His
work with the Hammond company, \citep{hammond1997}, on making a digital atlas which scales
properly, is another example of a visual problem.\footnote{The atlas
  also deals with the problem of non-overlapping labels, which is more
  a problem of statistical mechanics.} His favorite example was
a river which meanders around a city. On which side of the river will
the city be shown when you scale down the picture?
Another is his work on making bank notes which you can not photocopy
because they contain so many scales.
At least one of the many scales will be badly
sampled to the fixed pixel size of
any scanner. Another project was the question of how to guarantee
that the photograph of (say) a painting has the same colors as the
painting itself. Strictly speaking, this is not possible, since the
material is not the same. But Feigenbaum's idea was to actually
re-photograph a first printout, and then find the correct color mapping
by comparing the original photograph with the second one (For experts:
This problem does not have
a unique solution, and one gets better results
when there are many colors in the picture.)
Finally, he had a keen interest in photography, and
clearly, his love for what a camera (and the eye) can do infuses
his text.
Petzval's study (and design) of optics,
is another example showing his
fascination with the subject.
I mention all this because it shows a logical evolution
of his thinking about optics, vision and optical resolution which
finally led to his manuscript.

Feigenbaum worked for a very long time on this project, and, from the
oldest files he shared with me, I found  a program, called
anamorph.exe, dated July 12, 2006 (written in pascal), in which he
already programmed what he calls 3D-anamorphs (those are the
``correct'' anamorphs.)
For a modern version, see page \pageref{pag:html}.



\ssection{Evolution and Design}
\label{sec:evolutionanddesign}
\label{pag:ead}

The book had in principle an interesting section ``Evolution and
Design''
which had no text at Feigenbaum's death. However, we
  corresponded (and discussed) extensively about this subject in the
  summer of 2009. Feigenbaum insisted that the design of the fish eye
  is excellent and that the design of the land eye is
  awful. He also insisted that this means that evolution does not
  always find the best solution, even if there is a lot of time (since
  animals left water about 500 million years ago.) He also complained
  about the lack of good experimental measurements. The following
  paragraphs are taken from a letter of June, 18, 2009. I did not edit
  the text, because it shows a good example of Feigenbaum's thinking.
\medskip

\setlength{\leftskip}{1cm}
\setlength{\rightskip}{1cm}
``This matters.  The point is that I say I am showing that the simplest
well-designed optics for the task of air to water imaging already bests
evolution by almost an order of magnitude.  I have not gone here further in
actually optimizing, but evidently, it is hard to do worse than biology for
land eyes with any theoretical knowledge.  Now, this isn't quite true.  I've
explored replacing my one uniform spherical lens by a variable density one.
I can numerically assert that the uniform version is the optimal among them.
(I can't figure out how to prove this.)  In some ways, my simple design
appears to be the end of the line for this genre of optics with the
strongest refractor in front.  It is important then, without much more
elaborate designing, to be satisfied that we have already bested evolution.
This is what has bothered me, but that I now feel reasonably confident
about.  Too bad that someone hasn't built such a simple eye and
experimentally checked it out.

Anyway, the story for the fish eye is totally different.
  Here the design of
evolution is almost perfect:  It has optimized the design of optics within
its diffraction limits.  A careful fit to the crummy plot you've seen,
determines two radii at which they determined variations that they claim to
be optically significant.  However the better to fit to my analytic three
parameter family plus Gaussian bumps, shows the bumps are also where the
authors\footnote{In
the book, he analyzes the experiments of \citep{Kroger2001}. They show that the density in the
fish eye is \emph{not} constant. Feigenbaum then deduces what he
alludes to in his letter.} claim, and precisely where the 5th order caustic of the analytic
family has reached its final tangency to be within the self-determined
aperture.  Indeed the bumps are with the correct sign and size to now induce
a new cusp, turning it into a 7th order caustic.  The spill-over is within
the faster than exponential short diffractive tail of the caustic
diffraction, and so improves the lens from f/1.6 to at least f/1.5.  This is
impressive to discover merely by fitting.  Things are as right as could have
been ordered.

The reason is that the fish eye is better than a usual optimization problem.
By studying 5th order isotropic optics, it turns out that any mutation that
has the lens grow first with one index, and then uniformly with another is
win-win:  It either simultaneously improves both the brightness and
resolution, or simultaneously worsens both.  In the second case the fish is
dead in the water against predators and finding food, while the former is
highly favored.  Each further striation in radial density works the same
way.  This is why evolution made extravagantly good eyes for fish.  It
matters everything what the environment is, and what its ambient physics can
provide gratis.  Where physics is less forth-coming in its abundance, as in
the land eye case, evolution falls flat on its face, and simply constructs
engineering kluges.

This is a pr\'ecis of what I intend to say about what we have learned about
evolution from the comparative study of water and land eyes.  This is why I
need my land eye analysis to be impeccable.''

\setlength{\leftskip}{0cm}

\ssection{Using a hologram instead of the cylindrical tube}
\label{pag:holo}


I devised and have made a hologram--
made from a 3D anamorph--which shows the reflections from a
tube
in a quality 
similar to viewing an actual cylinder. The point is that holograms reproduce
the interferences of light waves and therefore, they are as good as
seeing a scene. Thus, they can capture what a normal camera can not
distinguish (as in \fref{fig:dots}) other than by changing focus.

The idea of making a hologram came already up in early 2015, in discussions
with Karl \xx{Knop}, a specialist in optics (among many other things). Feigenbaum was immediately interested in the
idea, and Knop was looking for somebody to make such a hologram, but
without success. I made a second attempt while editing the book.

Hologram making was very fashionable in the 1980s, but currently,
it is difficult to find professional 
hologram makers. The hologram
of the empire state building was made by Walter \xx{Spierings}, who is one 
of the leaders in this trade. Making the hologram is extremely
delicate, as the objects in question should not move by more that
1/20th of the wave length of light, about 50 nanometers ($0.5\cdot
10^{-7}$m). Hair has a diameter of about 50000 nanometers. This needs precise mounting of whatever is photographed,
but even the fluctuations of the density of air during exposition time
matter. Making a good hologram needs, among other parameters, a careful
control of the technical 
details of the optical setup. In the case at hand, laser speckle was
an issue, since the laser beam is reflected from the anamorph to the mirror.
The hologram is made onto photoresist 
(photolithography), whose surface consists of 
optical ridges at a submicron scale that will 
manipulate incoming lightwaves into 
reconstructing the 3D scene. From this master hologram
multiple copies can be made by embossing plastic 
with a metal shim made by galvanic means from the photoresist.

\begin{figure}[ht!]
  \CC{\includegraphics[width=\columnwidth]{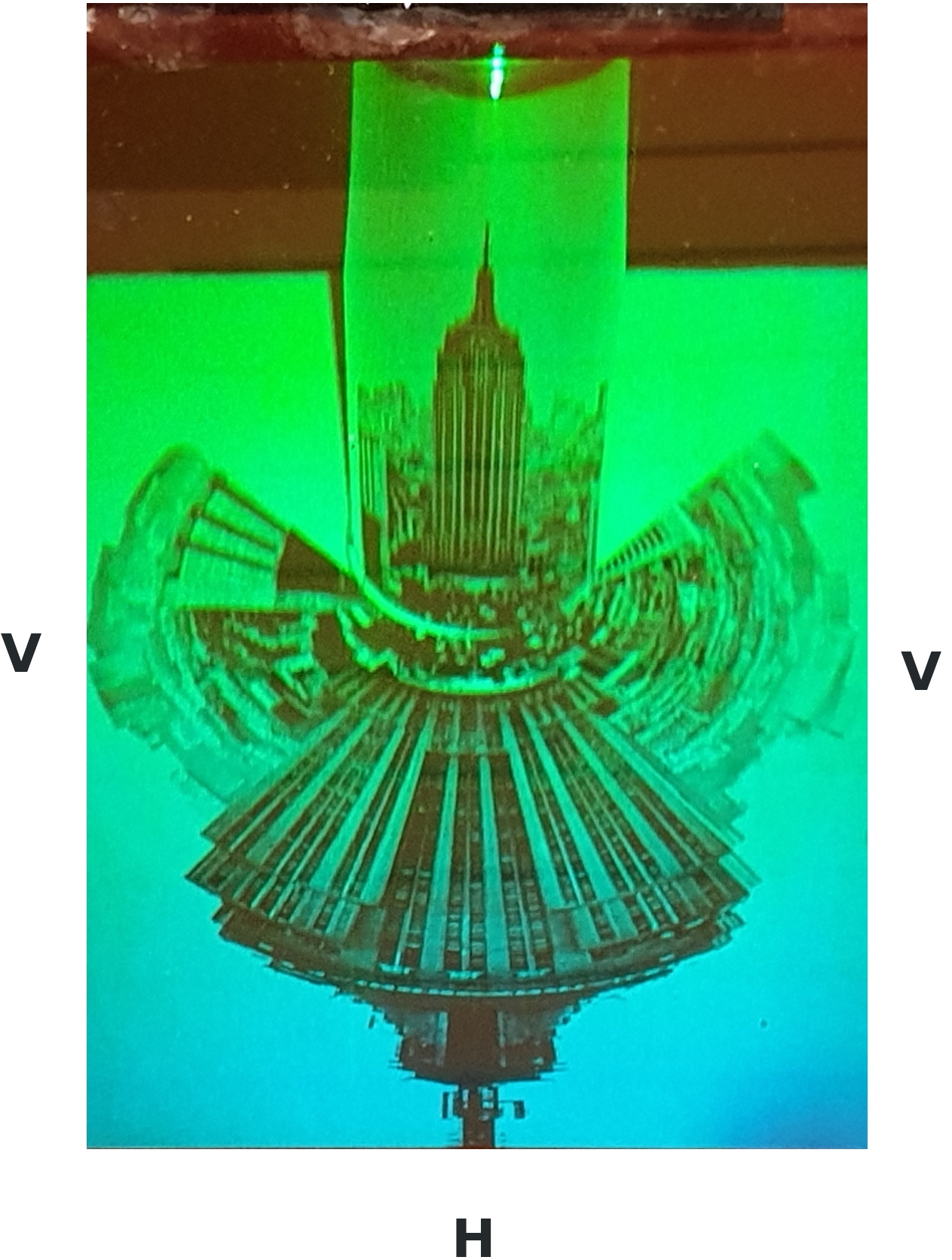}}
    \caption{Approximate rendering of the scene in the hologram. For
      best results, the observer should position the eye above the
      center of the hologram, looking down perpendicularly. If you
      stand on the {\textbf H} side you should see the empire state
      building in the tube, but standing on either of the {\textbf V}
      sides, you should see it in the same plane as the circular
      figure, as illustrated in \fref{fig:correct}. Note that the position of the image in the tube does not
    change.}\label{fig:holo}
    \end{figure}

\begin{figure}[ht!]
  \CC{\includegraphics[width=0.7\columnwidth]{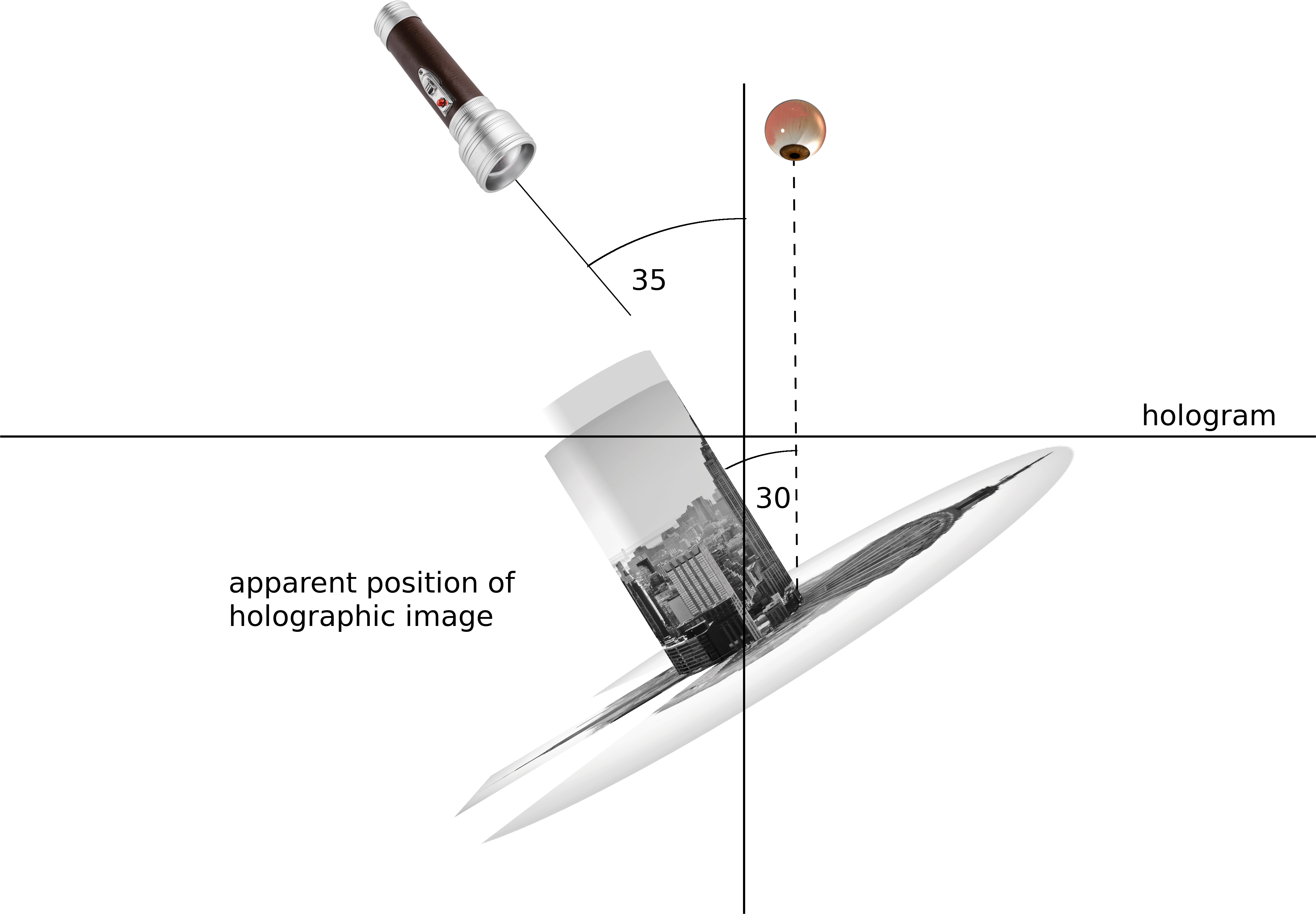}}
  \caption{A sketch of how the hologram will appear if the sheet is
    laid flat on a table. The illumination with a phone will work
    well. Results will be better if there is not too much stray light
  from other sources. Holding the light source at a larger distance is
  preferable. 
    The angle should be 35 degrees from the vertical (away from the
    top edge of the cardboard). Keep
  the light source fixed when viewing the hologram from the
  side. It is important that the light source is aligned with the axis
  of the hologram, and that the eye is really vertically above the
  center of the hologram. When you have the correct illumination, the
  hologram will appear uniformly green.}
  \label{fig:holo2}
  \end{figure}

\begin{figure}[ht!]
  {\includegraphics[width=0.64\columnwidth]{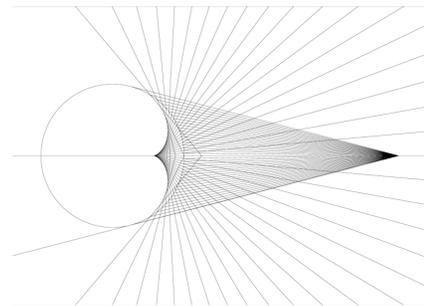}}
  \caption{The reflection of the rays emitted from a point as
    reflected from the circular cylinder. The caustic is inside the
    circle, and is the place where the rays are denser.}\label{fig:more1}
  
\end{figure}

The hologram has some advantages over viewing the scene with an actual
tube, apart from not needing to make a tube. In particular, instead of
turning the head one can just turn the 
hologram by 90 degrees. This avoids the neuro-ophthalmological signals
mentioned on page \pageref{pag:psycho}. The best results are obtained
if the hologram is illuminated from a fixed source, and then the
viewer walks around the table, as shown in \fref{fig:holo}. The
hologram is made in such a way that the best view is obtained when the
eye is perpendicularly above the center of the hologram as in \fref{fig:holo2}.

\ssection{Missing figures}

There were some pictures which had not made it into the manuscript. I
mention here again the one, \fref{fig:dots},
which ``proves'' the $\textbf{H}$ versus $\textbf{V}$
difference and which seems to me the most important:
It is a high quality photograph of a few of the
dots in \fref{fig:dots}.

Another image which should have appeared inside Feigenbaum's text is
the reflection of the rays in \fref{fig:more1}. Such drawings in 2
dimensions have a long history, but somehow, Feigenbaum's study in 3
dimensions seems to be absent in the literature.
\\

%


\nocite{jamescsouthall2011}
\nocite{hamilton1828}
\nocite{hamilton1830}
\nocite{hamilton1831}
\nocite{hamilton1833}
\nocite{hamilton1837}
\nocite{glaeser2000}
\nocite{synge1937}
\nocite{koretz1997}
\nocite{southall2010}
\nocite{thorne2017}
\nocite{sereno2008}
\nocite{smith1998}
\nocite{Morlet2002}
\nocite{campbell1966}
\nocite{guirao1999}
\nocite{bonds1974}
\nocite{Carkeet2003}
\nocite{adams1974}
\nocite{land1979}
\nocite{land2005}
\nocite{Nilsson1994}
\nocite{borghero2004}
\nocite{malmstrom2006}
\nocite{browne2007}
\nocite{zeiss2008}
\nocite{jagger1992}
\nocite{gordon200}
\nocite{abbe1878}
\nocite{okane2007}
\nocite{read2006}
\nocite{westheimer2005}
\nocite{westheimer1978}
\nocite{garding1995}
\nocite{qian1997}
\nocite{veltman1986}
\nocite{matthews2003}
\nocite{collewijn1991}
\nocite{1905s1}
\nocite{1905s2}
\nocite{1905s3}
\nocite{helmholtz1925}
\nocite{fernald}
\nocite{Bornwolf}
\nocite{stav1972}
\nocite{Shealy76}

\ssection{About references}
Feigenbaum's manuscript did not contain references. Fortunately, I could
find from my exchanges with him several files with relevant
literature. I decided to put them all as references at the end of the
this article,
whenever I could make out the source. The reader will find that
some of these references clearly address issues described in the planned
book. This collection also shows the eclectic interests of Feigenbaum.
He certainly was inspired by these references, but, knowing his way
of re-deriving everything by himself, we can assume that he did not
copy the results of the books and papers.
\medskip

\begin{acknowledgments}First of all, I am grateful to Mitchell
Feigenbaum to have shared, over the many years, his ideas, worries and
questions about the manuscript with me. I hope this outline gives justice
to his work. I also thank my many colleagues who commented on my
foreword, in particular, Michael Berry, David Campbell, Neil Dobbs, J\'er\'emie
Francfort, Gemunu Gunaratne, Karsten Kruse, Alberto
Morpurgo, Jacques Rougemont, and David Ruelle.
Some of their wishes could not be followed, because
I wanted to keep this summary of Feigenbaum's treatise as close as
possible to the spirit of the text.
No\'e Cuneo was instrumental in translating
Feigenbaum's ``Pascal'' program to html/javascript, so that readers can make
their own anamorphs.
I am grateful to Walter Spierings for his dedicated effort to make the
best possible hologram. I thank L\'eo L\'evy for help with making the index.
The cost of producing the holograms was kindly supported by the 
Fondation Schmidheiny, Geneva, and by the overhead of an ERC advanced 
grant, 290843 Bridges. 
\end{acknowledgments}
\bibliography{refs}

\end{document}